\documentclass[11pt]{article}
\usepackage{amsfonts}
\usepackage{bbm}
\usepackage{mathrsfs}
\usepackage{amssymb}

\def\hang{\hangindent\parindent}
\def\textindent#1{\indent\llap{#1\enspace}\ignorespaces}
\def\re{\par\hang\textindent}

\title{ Algebras Defined by Monic Gr\"obner Bases over Rings
\thanks{Project supported by
the National Natural Science Foundation of China (10571038,
10971044). \newline {\bf 2000 Mathematics Subject Classification:}
16W70; 16Z05. \newline {\bf Key words and phrases:}  Monic Gr\"obner
basis, Termination theorem, PBW $R$-basis, PBW isomorphism.}}

\author{Huishi Li\\
{\small Department of Applied Mathematics}\\
{\small College of Information Science and Technology}\\
{\small Hainan University}\\
{\small  Haikou 570228, China}}

\date{}

\begin{document}
\maketitle
\begin{center}
\begin{minipage}{120mm}
{\small {\bf Abstract.} Let $K\langle X\rangle =K\langle
X_1,...,X_n\rangle$ be the free algebra of $n$ generators over a
field $K$, and let $R\langle X\rangle =R\langle X_1,...,X_n\rangle$
be the free algebra of  $n$ generators over an arbitrary commutative
ring $R$. In this semi-expository paper, it is clarified that any
monic Gr\"obner basis in $K\langle X\rangle$ may give rise to a
monic Gr\"obner basis of the same type in $R\langle X\rangle$, and
vice versa. This fact turns out that many important $R$-algebras
have defining relations which form a monic Gr\"obner basis, and
consequently,  such $R$-algebras may be studied via a nice PBW
structure theory as that developed for quotient algebras of
$K\langle X\rangle$  in ([LWZ], [Li2, 3]). }
\end{minipage}\end{center}\par

\def\QED{\hfill{$\Box$}} \def\NZ{\mathbb{N}}
\def \r{\rightarrow}

\def\normalbaselines{\baselineskip 24pt\lineskip 4pt\lineskiplimit 4pt}
\def\mapdown#1{\llap{$\vcenter {\hbox {$\scriptstyle #1$}}$}
                                \Bigg\downarrow}
\def\mapdownr#1{\Bigg\downarrow\rlap{$\vcenter{\hbox
                                    {$\scriptstyle #1$}}$}}
\def\mapright#1#2{\smash{\mathop{\longrightarrow}\limits^{#1}_{#2}}}
\def\mapleft#1#2{\smash{\mathop{\longleftarrow}\limits^{#1}_{#2}}}
\def\mapup#1{\Bigg\uparrow\rlap{$\vcenter {\hbox  {$\scriptstyle #1$}}$}}
\def\mapupl#1{\llap{$\vcenter {\hbox {$\scriptstyle #1$}}$}
                                      \Bigg\uparrow}
\def\v5{\vskip .5truecm}
\def\T#1{\widetilde #1}
\def\OV#1{\overline {#1}}
\def\hang{\hangindent\parindent}
\def\textindent#1{\indent\llap{#1\enspace}\ignorespaces}
\def\item{\par\hang\textindent}
\message{<Paul Taylor's commutative diagrams, 20 July 1990>}
\newdimen\DiagramCellHeight\DiagramCellHeight3em 
\newdimen\DiagramCellWidth\DiagramCellWidth3em 
\newdimen\MapBreadth\MapBreadth.04em 
\newdimen\MapShortFall\MapShortFall.4em 
\newdimen\PileSpacing\PileSpacing1em 
\def\labelstyle{\ifincommdiag\textstyle\else\scriptstyle\fi}
\let\objectstyle\displaystyle


\def\rTo{\HorizontalMap\empty-\empty-\rhvee}
\def\lTo{\HorizontalMap\lhvee-\empty-\empty}
\def\dTo{\VerticalMap\empty|\empty|\dhvee}
\def\uTo{\VerticalMap\uhvee|\empty|\empty}
\let\uFrom\uTo\let\lFrom\lTo

\def\rArr{\HorizontalMap\empty-\empty-\rhla}
\def\lArr{\HorizontalMap\lhla-\empty-\empty}
\def\dArr{\VerticalMap\empty|\empty|\dhla}
\def\uArr{\VerticalMap\uhla|\empty|\empty}

\def\rDotsto{\HorizontalMap\empty\hfdot\hfdot\hfdot\rhvee}
\def\lDotsto{\HorizontalMap\lhvee\hfdot\hfdot\hfdot\empty}
\def\dDotsto{\VerticalMap\empty\vfdot\vfdot\vfdot\dhvee}
\def\uDotsto{\VerticalMap\uhvee\vfdot\vfdot\vfdot\empty}
\let\uDotsfrom\uDotsto\let\lDotsfrom\lDotsto

\def\rDashto{\HorizontalMap\empty\hfdash\hfdash\hfdash\rhvee}
\def\lDashto{\HorizontalMap\lhvee\hfdash\hfdash\hfdash\empty}
\def\dDashto{\VerticalMap\empty\vfdash\vfdash\vfdash\dhvee}
\def\uDashto{\VerticalMap\uhvee\vfdash\vfdash\vfdash\empty}
\let\uDashfrom\uDashto\let\lDashfrom\lDashto

\def\rImplies{\HorizontalMap==\empty=\Rightarrow}
\def\lImplies{\HorizontalMap\Leftarrow=\empty==}
\def\dImplies{\VerticalMap\|\|\empty\|\Downarrow}
\def\uImplies{\VerticalMap\Uparrow\|\empty\|\|}
\let\uImpliedby\uImplies\let\lImpliedby\lImplies

\def\rMapsto{\HorizontalMap\rtbar-\empty-\rhvee}
\def\lMapsto{\HorizontalMap\lhvee-\empty-\ltbar}
\def\dMapsto{\VerticalMap\dtbar|\empty|\dhvee}
\def\uMapsto{\VerticalMap\uhvee|\empty|\utbar}
\let\uMapsfrom\uMapsto\let\lMapsfrom\lMapsto

\def\rIntoA{\HorizontalMap\rthooka-\empty-\rhvee}
\def\rIntoB{\HorizontalMap\rthookb-\empty-\rhvee}
\def\lIntoA{\HorizontalMap\lhvee-\empty-\lthooka}
\def\lIntoB{\HorizontalMap\lhvee-\empty-\lthookb}
\def\dIntoA{\VerticalMap\dthooka|\empty|\dhvee}
\def\dIntoB{\VerticalMap\dthookb|\empty|\dhvee}
\def\uIntoA{\VerticalMap\uhvee|\empty|\uthooka}
\def\uIntoB{\VerticalMap\uhvee|\empty|\uthookb}
\let\uInfromA\uIntoA\let\uInfromB\uIntoB\let\lInfromA\lIntoA\let\lInfromB
\lIntoB\let\rInto\rIntoA\let\lInto\lIntoA\let\dInto\dIntoB\let\uInto\uIntoA

\def\rEmbed{\HorizontalMap\gt-\empty-\rhvee}
\def\lEmbed{\HorizontalMap\lhvee-\empty-\lt}
\def\dEmbed{\VerticalMap\vee|\empty|\dhvee}
\def\uEmbed{\VerticalMap\uhvee|\empty|\wedge}

\def\rProject{\HorizontalMap\empty-\empty-\triangleright}
\def\lProject{\HorizontalMap\triangleleft-\empty-\empty}
\def\uProject{\VerticalMap\triangleup|\empty|\empty}
\def\dProject{\VerticalMap\empty|\empty|\littletriangledown}

\def\rOnto{\HorizontalMap\empty-\empty-\twoheadrightarrow}
\def\lOnto{\HorizontalMap\twoheadleftarrow-\empty-\empty}
\def\dOnto{\VerticalMap\empty|\empty|\twoheaddownarrow}
\def\uOnto{\VerticalMap\twoheaduparrow|\empty|\empty}
\let\lOnfrom\lOnto\let\uOnfrom\uOnto

\def\hEq{\HorizontalMap==\empty==}
\def\vEq{\VerticalMap\|\|\empty\|\|}
\let\rEq\hEq\let\lEq\hEq\let\uEq\vEq\let\dEq\vEq

\def\hLine{\HorizontalMap\empty-\empty-\empty}
\def\vLine{\VerticalMap\empty|\empty|\empty}
\let\rLine\hLine\let\lLine\hLine\let\uLine\vLine\let\dLine\vLine

\def\hDots{\HorizontalMap\empty\hfdot\hfdot\hfdot\empty}
\def\vDots{\VerticalMap\empty\vfdot\vfdot\vfdot\empty}
\let\rDots\hDots\let\lDots\hDots\let\uDots\vDots\let\dDots\vDots

\def\hDashes{\HorizontalMap\empty\hfdash\hfdash\hfdash\empty}
\def\vDashes{\VerticalMap\empty\vfdash\vfdash\vfdash\empty}
\let\rDashes\hDashes\let\lDashes\hDashes\let\uDashes\vDashes\let\dDashes
\vDashes

\def\rPto{\HorizontalMap\empty-\empty-\rightharpoonup}
\def\lPto{\HorizontalMap\leftharpoonup-\empty-\empty}
\def\uPto{\VerticalMap\upharpoonright|\empty|\empty}
\def\dPto{\VerticalMap\empty|\empty|\downharpoonright}
\let\lPfrom\lPto\let\uPfrom\uPto

\def\NW{\NorthWest\DiagonalMap{\lah111}{\laf100}{}{\laf100}{}(2,2)}
\def\NE{\NorthEast\DiagonalMap{\lah22}{\laf0}{}{\laf0}{}(2,2)}
\def\SW{\SouthWest\DiagonalMap{}{\laf0}{}{\laf0}{\lah11}(2,2)}
\def\SE{\SouthEast\DiagonalMap{}{\laf100}{}{\laf100}{\lah122}(2,2)}

\def\nNW{\NorthWest\DiagonalMap{\lah135}{\laf112}{}{\laf112}{}(2,3)}
\def\nNE{\NorthEast\DiagonalMap{\lah36}{\laf12}{}{\laf12}{}(2,3)}
\def\sSW{\SouthWest\DiagonalMap{}{\laf12}{}{\laf12}{\lah35}(2,3)}
\def\sSE{\SouthEast\DiagonalMap{}{\laf112}{}{\laf112}{\lah136}(2,3)}

\def\wNW{\NorthWest\DiagonalMap{\lah153}{\laf121}{}{\laf121}{}(3,2)}
\def\eNE{\NorthEast\DiagonalMap{\lah63}{\laf21}{}{\laf21}{}(3,2)}
\def\wSW{\SouthWest\DiagonalMap{}{\laf21}{}{\laf21}{\lah53}(3,2)}
\def\eSE{\SouthEast\DiagonalMap{}{\laf121}{}{\laf121}{\lah163}(3,2)}

\def\NNW{\NorthWest\DiagonalMap{\lah113}{\laf101}{}{\laf101}{}(2,4)}
\def\NNE{\NorthEast\DiagonalMap{\lah25}{\laf01}{}{\laf01}{}(2,4)}
\def\SSW{\SouthWest\DiagonalMap{}{\laf01}{}{\laf01}{\lah13}(2,4)}
\def\SSE{\SouthEast\DiagonalMap{}{\laf101}{}{\laf101}{\lah125}(2,4)}

\def\WNW{\NorthWest\DiagonalMap{\lah131}{\laf110}{}{\laf110}{}(4,2)}
\def\ENE{\NorthEast\DiagonalMap{\lah52}{\laf10}{}{\laf10}{}(4,2)}
\def\WSW{\SouthWest\DiagonalMap{}{\laf10}{}{\laf10}{\lah31}(4,2)}
\def\ESE{\SouthEast\DiagonalMap{}{\laf110}{}{\laf110}{\lah152}(4,2)}

\def\NNNW{\NorthWest\DiagonalMap{\lah115}{\laf102}{}{\laf102}{}(2,6)}
\def\NNNE{\NorthEast\DiagonalMap{\lah16}{\laf02}{}{\laf02}{}(2,6)}
\def\SSSW{\SouthWest\DiagonalMap{}{\laf02}{}{\laf02}{\lah15}(2,6)}
\def\SSSE{\SouthEast\DiagonalMap{}{\laf102}{}{\laf102}{\lah116}(2,6)}

\def\WWNW{\NorthWest\DiagonalMap{\lah151}{\laf120}{}{\laf120}{}(6,2)}
\def\EENE{\NorthEast\DiagonalMap{\lah61}{\laf20}{}{\laf20}{}(6,2)}
\def\WWSW{\SouthWest\DiagonalMap{}{\laf20}{}{\laf20}{\lah51}(6,2)}
\def\EESE{\SouthEast\DiagonalMap{}{\laf120}{}{\laf120}{\lah161}(6,2)}


\font\tenln=line10

\mathchardef\lt="313C \mathchardef\gt="313E

\def\rhvee{\mkern-10mu\gt}
\def\lhvee{\lt\mkern-10mu}
\def\dhvee{\vbox\tozpt{\vss\hbox{$\vee$}\kern0pt}}
\def\uhvee{\vbox\tozpt{\hbox{$\wedge$}\vss}}
\def\rhcvee{\mkern-10mu\succ}
\def\lhcvee{\prec\mkern-10mu}
\def\dhcvee{\vbox\tozpt{\vss\hbox{$\curlyvee$}\kern0pt}}
\def\uhcvee{\vbox\tozpt{\hbox{$\curlywedge$}\vss}}
\def\rhvvee{\mkern-10mu\gg}
\def\lhvvee{\ll\mkern-10mu}
\def\dhvvee{\vbox\tozpt{\vss\hbox{$\vee$}\kern-.6ex\hbox{$\vee$}\kern0pt}}
\def\uhvvee{\vbox\tozpt{\hbox{$\wedge$}\kern-.6ex\hbox{$\wedge$}\vss}}
\def\twoheaddownarrow{\rlap{$\downarrow$}\raise-.5ex\hbox{$\downarrow$}}
\def\twoheaduparrow{\rlap{$\uparrow$}\raise.5ex\hbox{$\uparrow$}}
\def\triangleup{{\scriptscriptstyle\bigtriangleup}}
\def\littletriangledown{{\scriptscriptstyle\triangledown}}
\def\rhla{\vbox\tozpt{\vss\hbox\tozpt{\hss\tenln\char'55}\kern\axisheight}}
\def\lhla{\vbox\tozpt{\vss\hbox\tozpt{\tenln\char'33\hss}\kern\axisheight}}
\def\dhla{\vbox\tozpt{\vss\hbox\tozpt{\tenln\char'77\hss}}}
\def\uhla{\vbox\tozpt{\hbox\tozpt{\tenln\char'66\hss}\vss}}
\def\htdot{\mkern3.15mu\cdot\mkern3.15mu}
\def\vtdot{\vbox to 1.46ex{\vss\hbox{$\cdot$}}}
\def\utbar{\vrule height 0.093ex depth0pt width 0.4em} \let\dtbar\utbar
\def\rtbar{\mkern1.5mu\vrule height 1.1ex depth.06ex width .04em\mkern1.5mu}%
\let\ltbar\rtbar
\def\rthooka{\raise.603ex\hbox{$\scriptscriptstyle\subset$}}
\def\lthooka{\raise.603ex\hbox{$\scriptscriptstyle\supset$}}
\def\rthookb{\raise-.022ex\hbox{$\scriptscriptstyle\subset$}}
\def\lthookb{\raise-.022ex\hbox{$\scriptscriptstyle\supset$}}
\def\dthookb{\hbox{$\scriptscriptstyle\cap$}\mkern5.5mu}
\def\uthookb{\hbox{$\scriptscriptstyle\cup$}\mkern4.5mu}
\def\dthooka{\mkern6mu\hbox{$\scriptscriptstyle\cap$}}
\def\uthooka{\mkern6mu\hbox{$\scriptscriptstyle\cup$}}
\def\hfdot{\mkern3.15mu\cdot\mkern3.15mu}
\def\vfdot{\vbox to 1.46ex{\vss\hbox{$\cdot$}}}
\def\vfdashstrut{\vrule width0pt height1.3ex depth0.7ex}
\def\vfthedash{\vrule width\MapBreadth height0.6ex depth 0pt}
\def\hfthedash{\vrule\horizhtdp width 0.26em}
\def\hfdash{\mkern5.5mu\hfthedash\mkern5.5mu}
\def\vfdash{\vfdashstrut\vfthedash}

\def\nwhTO{\nwarrow\mkern-1mu}
\def\nehTO{\mkern-.1mu\nearrow}
\def\sehTO{\searrow\mkern-.02mu}
\def\swhTO{\mkern-.8mu\swarrow}

\def\SEpbk{\rlap{\smash{\kern0.1em \vrule depth 2.67ex height -2.55ex width 0%
.9em \vrule height -0.46ex depth 2.67ex width .05em }}}
\def\SWpbk{\llap{\smash{\vrule height -0.46ex depth 2.67ex width .05em \vrule
depth 2.67ex height -2.55ex width .9em \kern0.1em }}}
\def\NEpbk{\rlap{\smash{\kern0.1em \vrule depth -3.48ex height 3.67ex width 0%
.95em \vrule height 3.67ex depth -1.39ex width .05em }}}
\def\NWpbk{\llap{\smash{\vrule height 3.6ex depth -1.39ex width .05em \vrule
depth -3.48ex height 3.67ex width .95em \kern0.1em }}}


\newcount\cdna\newcount\cdnb\newcount\cdnc\newcount\cdnd\cdna=\catcode`\@%
\catcode`\@=11 \let\then\relax\def\loopa#1\repeat{\def\bodya{#1}\iteratea}%
\def\iteratea{\bodya\let\next\iteratea\else\let\next\relax\fi\next}\def\loopb
#1\repeat{\def\bodyb{#1}\iterateb}\def\iterateb{\bodyb\let\next\iterateb\else
\let\next\relax\fi\next} \def\mapctxterr{\message{commutative diagram: map
context error}}\def\mapclasherr{\message{commutative diagram: clashing maps}}%
\def\ObsDim#1{\expandafter\message{! diagrams Warning: Dimension \string#1 is
obsolete
(ignored)}\global\let#1\ObsDimq\ObsDimq}\def\ObsDimq{\dimen@=}\def
\HorizontalMapLength{\ObsDim\HorizontalMapLength}\def\VerticalMapHeight{%
\ObsDim\VerticalMapHeight}\def\VerticalMapDepth{\ObsDim\VerticalMapDepth}\def
\VerticalMapExtraHeight{\ObsDim\VerticalMapExtraHeight}\def
\VerticalMapExtraDepth{\ObsDim\VerticalMapExtraDepth}\def\ObsCount#1{%
\expandafter\message{! diagrams Warning: Count \string#1 is obsolete (ignored%
)}\global\let#1\ObsCountq\ObsCountq}\def\ObsCountq{\count@=}\def
\DiagonalLineSegments{\ObsCount\DiagonalLineSegments}\def\tozpt{to\z@}\def
\sethorizhtdp{\dimen8=\axisheight\dimen9=\MapBreadth\advance\dimen8.5\dimen9%
\advance\dimen9-\dimen8}\def\horizhtdp{height\dimen8 depth\dimen9
}\def \axisheight{\fontdimen22\the\textfont2 }\countdef\boxc@unt=14

\def\bombparameters{\hsize\z@\rightskip\z@ plus1fil minus\maxdimen
\parfillskip\z@\linepenalty9000 \looseness0 \hfuzz\maxdimen\hbadness10000
\clubpenalty0 \widowpenalty0 \displaywidowpenalty0
\interlinepenalty0 \predisplaypenalty0 \postdisplaypenalty0
\interdisplaylinepenalty0
\interfootnotelinepenalty0 \floatingpenalty0 \brokenpenalty0 \everypar{}%
\leftskip\z@\parskip\z@\parindent\z@\pretolerance10000
\tolerance10000 \hyphenpenalty10000 \exhyphenpenalty10000
\binoppenalty10000 \relpenalty10000 \adjdemerits0
\doublehyphendemerits0 \finalhyphendemerits0 \prevdepth\z@}\def
\startbombverticallist{\hbox{}\penalty1\nointerlineskip}

\def\pushh#1\to#2{\setbox#2=\hbox{\box#1\unhbox#2}}\def\pusht#1\to#2{\setbox#%
2=\hbox{\unhbox#2\box#1}}

\newif\ifallowhorizmap\allowhorizmaptrue\newif\ifallowvertmap
\allowvertmapfalse\newif\ifincommdiag\incommdiagfalse

\def\diagram{\hbox\bgroup$\vcenter\bgroup\startbombverticallist
\incommdiagtrue\baselineskip\DiagramCellHeight\lineskip\z@\lineskiplimit\z@
\mathsurround\z@\tabskip\z@\let\\\diagcr\allowhorizmaptrue\allowvertmaptrue
\halign\bgroup\lcdtempl##\rcdtempl&&\lcdtempl##\rcdtempl\cr}\def\enddiagram{%
\crcr\egroup\reformatmatrix\egroup$\egroup}\def\commdiag#1{{\diagram#1%
\enddiagram}}

\def\lcdtempl{\futurelet\thefirsttoken\dolcdtempl}\newif\ifemptycell\def
\dolcdtempl{\ifx\thefirsttoken\rcdtempl\then\hskip1sp plus 1fil
\emptycelltrue
\else\hfil$\emptycellfalse\objectstyle\fi}\def\rcdtempl{\ifemptycell\else$%
\hfil\fi}\def\diagcr{\cr} \def\across#1{\span\omit\mscount=#1
\loop\ifnum \mscount>2
\spAn\repeat\ignorespaces}\def\spAn{\relax\span\omit\advance
\mscount by -1}

\def\CellSize{\afterassignment\cdhttowd\DiagramCellHeight}\def\cdhttowd{%
\DiagramCellWidth\DiagramCellHeight}\def\MapsAbut{\MapShortFall\z@}

\newcount\cdvdl\newcount\cdvdr\newcount\cdvd\newcount\cdbfb\newcount\cdbfr
\newcount\cdbfl\newcount\cdvdr\newcount\cdvdl\newcount\cdvd

\def\reformatmatrix{\bombparameters\cdvdl=\insc@unt\cdvdr=\cdvdl\cdbfb=%
\boxc@unt\advance\cdbfb1
\cdbfr=\cdbfb\setbox1=\vbox{}\dimen2=\z@\loop\setbox
0=\lastbox\ifhbox0 \dimen1=\lastskip\unskip\dimen5=\ht0
\advance\dimen5 \dimen
1 \dimen4=\dp0 \penalty1 \reformatrow\unpenalty\ht4=\dimen5 \dp4=\dimen4 \ht3%
\z@\dp3\z@\setbox1=\vbox{\box4 \nointerlineskip\box3 \nointerlineskip\unvbox1%
}\dimen2=\dimen1 \repeat\unvbox1}

\newif\ifcontinuerow

\def\reformatrow{\cdbfl=\cdbfr\noindent\unhbox0 \loopa\unskip\setbox\cdbfl=%
\lastbox\ifhbox\cdbfl\advance\cdbfl1\repeat\par\unskip\dimen6=2%
\DiagramCellWidth\dimen7=-\DiagramCellWidth\setbox3=\hbox{}\setbox4=\hbox{}%
\setbox7=\box\voidb@x\cdvd=\cdvdl\continuerowtrue\loopa\advance\cdvd-1
\adjustcells\ifcontinuerow\advance\dimen6\wd\cdbfl\cdda=.5\dimen6
\ifdim\cdda
<\DiagramCellWidth\then\dimen6\DiagramCellWidth\advance\dimen6-\cdda
\nopendvert\cdda\DiagramCellWidth\fi\advance\dimen7\cdda\dimen6=\wd\cdbfl
\reformatcell\advance\cdbfl-1 \repeat\advance\dimen7.5\dimen6
\outHarrow} \def
\adjustcells{\ifnum\cdbfr>\cdbfl\then\ifnum\cdvdr>\cdvd\then\continuerowfalse
\else\setbox\cdbfl=\hbox
to\wd\cdvd{\lcdtempl\VonH{}\rcdtempl}\fi\else\ifnum
\cdvdr>\cdvd\then\advance\cdvdr-1
\setbox\cdvd=\vbox{}\wd\cdvd=\wd\cdbfl\dp \cdvd=\dp1 \fi\fi}

\def\reformatcell{\sethorizhtdp\noindent\unhbox\cdbfl\skip0=\lastskip\unskip
\par\ifcase\prevgraf\reformatempty\or\reformatobject\else\reformatcomplex\fi
\unskip}\def\reformatobject{\setbox6=\lastbox\unskip\vadjdon6\outVarrow
\setbox6=\hbox{\unhbox6}\advance\dimen7-.5\wd6
\outHarrow\dimen7=-.5\wd6 \pusht6\to4}\newcount\globnum

\def\reformatcomplex{\setbox6=\lastbox\unskip\setbox9=\lastbox\unskip\setbox9%
=\hbox{\unhbox9
\skip0=\lastskip\unskip\global\globnum\lastpenalty\hskip\skip 0
}\advance\globnum9999
\ifcase\globnum\reformathoriz\or\reformatpile\or
\reformatHonV\or\reformatVonH\or\reformatvert\or\reformatHmeetV\fi}

\def\reformatempty{\vpassdon\ifdim\skip0>\z@\then\hpassdon\else\ifvoid2 \then
\else\advance\dimen7-.5\dimen0
\cdda=\wd2\advance\cdda.5\dimen0\wd2=\cdda\fi
\fi}\def\VonH{\doVonH6}\def\HonV{\doVonH7}\def\HmeetV{\MapBreadth-2%
\MapShortFall\doVonH4}\def\doVonH#1{\cdna-999#1\futurelet\thenexttoken
\dooVonH}\def\dooVonH{\let\next\relax\sethorizhtdp\ifallowhorizmap
\ifallowvertmap\then\ifx\thenexttoken[\then\let\next\VonHstrut\else
\sethorizhtdp\dimen0\MapBreadth\let\next\VonHnostrut\fi\else\mapctxterr\fi
\else\mapctxterr\fi\next}\def\VonHstrut[#1]{\setbox0=\hbox{$#1$}\dimen0\wd0%
\dimen8\ht0\dimen9\dp0 \VonHnostrut}\def\VonHnostrut{\setbox0=\hbox{}\ht0=%
\dimen8\dp0=\dimen9\wd0=.5\dimen0 \copy0\penalty\cdna\box0
\allowhorizmapfalse
\allowvertmapfalse}\def\reformatHonV{\hpassdon\doreformatHonV}\def
\reformatHmeetV{\dimen@=\wd9 \advance\dimen7-\wd9 \outHarrow\setbox6=\hbox{%
\unhbox6}\dimen7-\wd6 \advance\dimen@\wd6 \setbox6=\hbox to\dimen@{\hss}%
\pusht6\to4\doreformatHonV}\def\doreformatHonV{\setbox9=\hbox{\unhbox9
\unskip
\unpenalty\global\setbox\globbox=\lastbox}\vadjdon\globbox\outVarrow}\def
\reformatVonH{\vpassdon\advance\dimen7-\wd9 \outHarrow\setbox6=\hbox{\unhbox6%
}\dimen7=-\wd6 \setbox6=\hbox{\kern\wd9 \kern\wd6}\pusht6\to4}\def\hpassdon{}%
\def\vpassdon{\dimen@=\dp\cdvd\advance\dimen@\dimen4 \advance\dimen@\dimen5
\dp\cdvd=\dimen@\nopendvert}\def\vadjdon#1{\dimen8=\ht#1
\dimen9=\dp#1 }

\def\HorizontalMap#1#2#3#4#5{\sethorizhtdp\setbox1=\makeharrowpart{#1}\def
\arrowfillera{#2}\def\arrowfillerb{#4}\setbox5=\makeharrowpart{#5}\ifx
\arrowfillera\justhorizline\then\def\arra{\hrule\horizhtdp}\def\kea{\kern-0.%
01em}\let\arrstruthtdp\horizhtdp\else\def\kea{\kern-0.15em}\setbox2=\hbox{%
\kea${\arrowfillera}$\kea}\def\arra{\copy2}\def\arrstruthtdp{height\ht2 depth%
\dp2
}\fi\ifx\arrowfillerb\justhorizline\then\def\arrb{\hrule\horizhtdp}\def
\keb{kern-0.01em}\ifx\arrowfillera\empty\then\let\arrstruthtdp\horizhtdp\fi
\else\def\keb{\kern-0.15em}\setbox4=\hbox{\keb${\arrowfillerb}$\keb}\def\arrb
{\copy4}\ifx\arrowfilera\empty\then\def\arrstruthtdp{height\ht4 depth\dp4 }%
\fi\fi\setbox3=\makeharrowpart{{#3}\vrule width\z@\arrstruthtdp}%
\ifallowhorizmap\then\let\execmap\execHorizontalMap\else\let\execmap
\mapctxterr\fi\allowhorizmapfalse\gettwoargs}\def\makeharrowpart#1{\hbox{%
\mathsurround\z@\edef\next{#1}\ifx\next\empty\else$\mkern-1.5mu{\next}\mkern-%
1.5mu$\fi}}\def\justhorizline{-}

\def\execHorizontalMap{\dimen0=\wd6 \ifdim\dimen0<\wd7\then\dimen0=\wd7\fi
\dimen3=\wd3 \ifdim\dimen0<2em\then\dimen0=2em\fi\skip2=.5\dimen0
\ifincommdiag plus 1fill\fi minus\z@\advance\skip2-.5\dimen3
\skip4=\skip2 \advance\skip2-\wd1 \advance\skip4-\wd5
\kern\MapShortFall\box1 \xleaders
\arra\hskip\skip2 \vbox{\lineskiplimit\maxdimen\lineskip.5ex \ifhbox6 \hbox to%
\dimen3 {\hss\box6\hss}\fi\vtop{\box3 \ifhbox7 \hbox to\dimen3
{\hss\box7\hss
}\fi}}\ifincommdiag\kern-.5\dimen3\penalty-9999\null\kern.5\dimen3\fi
\xleaders\arrb\hskip\skip4 \box5 \kern\MapShortFall}

\def\reformathoriz{\vadjdon6\outVarrow\ifvoid7\else\mapclasherr\fi\setbox2=%
\box9 \wd2=\dimen7 \dimen7=\z@\setbox7=\box6 }

\def\resetharrowpart#1#2{\ifvoid#1\then\ifdim#2=\z@\else\setbox4=\hbox{%
\unhbox4\kern#2}\fi\else\ifhbox#1\then\setbox#1=\hbox
to#2{\unhbox#1}\else
\widenpile#1\fi\pusht#1\to4\fi}\def\outHarrow{\resetharrowpart2{\wd2}\pusht2%
\to4\resetharrowpart7{\dimen7}\pusht7\to4\dimen7=\z@}

\def\pile#1{{\incommdiagtrue\let\pile\innerpile\allowvertmapfalse
\allowhorizmaptrue\baselineskip.5\PileSpacing\lineskip\z@\lineskiplimit\z@
\mathsurround\z@\tabskip\z@\let\\\pilecr\vcenter{\halign{\hfil$##$\hfil\cr#1
\crcr}}}\ifincommdiag\then\ifallowhorizmap\then\penalty-9998
\allowvertmapfalse\allowhorizmapfalse\else\mapctxterr\fi\fi}\def\pilecr{\cr}%
\def\innerpile#1{\noalign{\halign{\hfil$##$\hfil\cr#1 \crcr}}}

\def\reformatpile{\vadjdon9\outVarrow\ifvoid7\else\mapclasherr\fi\penalty1
\setbox9=\hbox{\unhbox9 \unskip\unpenalty\setbox9=\lastbox\unhbox9
\global
\setbox\globbox=\lastbox}\unvbox\globbox\setbox9=\vbox{}\setbox7=\vbox{}%
\loopb\setbox6=\lastbox\ifhbox6
\skip3=\lastskip\unskip\splitpilerow\repeat
\unpenalty\setbox9=\hbox{$\vcenter{\unvbox9}$}\setbox2=\box9
\dimen7=\z@}\def
\pilestrut{\vrule height\dimen0 depth\dimen3 width\z@}\def\splitpilerow{%
\dimen0=\ht6 \dimen3=\dp6
\noindent\unhbox6\unskip\setbox6=\lastbox\unskip
\unhbox6\par\setbox6=\lastbox\unskip\ifcase\prevgraf\or\setbox6=\hbox\tozpt{%
\hss\unhbox6\hss}\ht6=\dimen0 \dp6=\dimen3
\setbox9=\vbox{\vskip\skip3 \hbox
to\dimen7{\hfil\box6}\nointerlineskip\unvbox9}\setbox7=\vbox{\vskip\skip3
\hbox{\pilestrut\hfil}\nointerlineskip\unvbox7}\or\setbox7=\vbox{\vskip\skip3
\hbox{\pilestrut\unhbox6}\nointerlineskip\unvbox7}\setbox6=\lastbox\unskip
\setbox9=\vbox{\vskip\skip3 \hbox to\dimen7{\pilestrut\unhbox6}%
\nointerlineskip\unvbox9}\fi\unskip}

\def\widenpile#1{\setbox#1=\hbox{$\vcenter{\unvbox#1 \setbox8=\vbox{}\loopb
\setbox9=\lastbox\ifhbox9
\skip3=\lastskip\unskip\setbox8=\vbox{\vskip\skip3 \hbox
to\dimen7{\unhbox9}\nointerlineskip\unvbox8}\repeat\unvbox8 }$}}

\def\justverticalline{|}\def\makevarrowpart#1{\hbox to\MapBreadth{\hss$\kern
\MapBreadth{#1}$\hss}}\def\VerticalMap#1#2#3#4#5{\setbox1=\makevarrowpart{#1}%
\def\arrowfillera{#2}\setbox3=\makevarrowpart{#3}\def\arrowfillerb{#4}\setbox
5=\makevarrowpart{#5}\ifx\arrowfillera\justverticalline\then\def\arra{\vrule
width\MapBreadth}\def\kea{\kern-0.05ex}\else\def\kea{\kern-0.35ex}\setbox2=%
\vbox{\kea\makevarrowpart\arrowfillera\kea}\def\arra{\copy2}\fi\ifx
\arrowfillerb\justverticalline\then\def\arrb{\vrule
width\MapBreadth}\def\keb
{\kern-0.05ex}\else\def\keb{\kern-0.35ex}\setbox4=\vbox{\keb\makevarrowpart
\arrowfillerb\keb}\def\arrb{\copy4}\fi\ifallowvertmap\then\let\execmap
\execVerticalMap\else\let\execmap\mapctxterr\fi\allowhorizmapfalse\gettwoargs
}

\def\execVerticalMap{\setbox3=\makevarrowpart{\box3}\setbox0=\hbox{}\ht0=\ht3
\dp0\z@\ht3\z@\box6 \setbox8=\vtop spread2ex{\offinterlineskip\box3
\xleaders
\arrb\vfill\box5 \kern\MapShortFall}\dp8=\z@\box8 \kern-\MapBreadth\setbox8=%
\vbox spread2ex{\offinterlineskip\kern\MapShortFall\box1
\xleaders\arra\vfill \box0}\ht8=\z@\box8
\ifincommdiag\then\kern-.5\MapBreadth\penalty-9995 \null
\kern.5\MapBreadth\fi\box7\hfil}

\newcount\colno\newdimen\cdda\newbox\globbox\def\reformatvert{\setbox6=\hbox{%
\unhbox6}\cdda=\wd6 \dimen3=\dp\cdvd\advance\dimen3\dimen4
\setbox\cdvd=\hbox {}\colno=\prevgraf\advance\colno-2
\loopb\setbox9=\hbox{\unhbox9 \unskip
\unpenalty\dimen7=\lastkern\unkern\global\setbox\globbox=\lastbox\advance
\dimen7\wd\globbox\advance\dimen7\lastkern\unkern\setbox9=\lastbox\vtop to%
\dimen3{\unvbox9}\kern\dimen7 }\ifnum\colno>0
\ifdim\wd9<\PileSpacing\then \setbox9=\hbox
to\PileSpacing{\unhbox9}\fi\dimen0=\wd9 \advance\dimen0-\wd
\globbox\setbox\cdvd=\hbox{\kern\dimen0 \box\globbox\unhbox\cdvd}\pushh9\to6%
\advance\colno-1
\setbox9=\lastbox\unskip\repeat\advance\dimen7-.5\wd6
\advance\dimen7.5\cdda\advance\dimen7-\wd9 \outHarrow\dimen7=-.5\wd6
\advance
\dimen7-.5\cdda\pusht9\to4\pusht6\to4\nopendvert\dimen@=\dimen6\advance
\dimen@-\wd\cdvd\advance\dimen@-\wd\globbox\divide\dimen@2
\setbox\cdvd=\hbox
{\kern\dimen@\box\globbox\unhbox\cdvd\kern\dimen@}\dimen8=\dp\cdvd\advance
\dimen8\dimen5 \dp\cdvd=\dimen8 \ht\cdvd=\z@}

\def\outVarrow{\ifhbox\cdvd\then\deepenbox\cdvd\pusht\cdvd\to3\else
\nopendvert\fi\dimen3=\dimen5 \advance\dimen3-\dimen8 \setbox\cdvd=\vbox{%
\vfil}\dp\cdvd=\dimen3} \def\nopendvert{\setbox3=\hbox{\unhbox3\kern\dimen6}}%
\def\deepenbox\cdvd{\setbox\cdvd=\hbox{\dimen3=\dimen4 \advance\dimen3-\dimen
9 \setbox6=\hbox{}\ht6=\dimen3 \dp6=-\dimen3 \dimen0=\dp\cdvd\advance\dimen0%
\dimen3 \unhbox\cdvd\dimen3=\lastkern\unkern\setbox8=\hbox{\kern\dimen3}%
\loopb\setbox9=\lastbox\ifvbox9 \setbox9=\vtop to\dimen0{\copy6
\nointerlineskip\unvbox9 }\dimen3=\lastkern\unkern\setbox8=\hbox{\kern\dimen3%
\box9\unhbox8}\repeat\unhbox8 }}

\newif\ifPositiveGradient\PositiveGradienttrue\newif\ifClimbing\Climbingtrue
\newcount\DiagonalChoice\DiagonalChoice1 \newcount\lineno\newcount\rowno
\newcount\charno\def\laf{\afterassignment\xlaf\charno='}\def\xlaf{\hbox{%
\tenln\char\charno}}\def\lah{\afterassignment\xlah\charno='}\def\xlah{\hbox{%
\tenln\char\charno}}\def\makedarrowpart#1{\hbox{\mathsurround\z@${#1}$}}\def
\lad{\afterassignment\xlad\charno='}\def\xlad{\setbox2=\xlaf\setbox0=\hbox to%
.5\wd2{$\hss\ldot\hss$}\ht0=.25\ht2 \dp0=\ht0 \hbox{\mv-\ht0\copy0 \mv\ht0%
\box0}}

\def\DiagonalMap#1#2#3#4#5{\ifPositiveGradient\then\let\mv\raise\else\let\mv
\lower\fi\setbox2=\makedarrowpart{#2}\setbox1=\makedarrowpart{#1}\setbox4=%
\makedarrowpart{#4}\setbox5=\makedarrowpart{#5}\setbox3=\makedarrowpart{#3}%
\let\execmap\execDiagonalLine\gettwoargs}

\def\makeline#1(#2,#3;#4){\hbox{\dimen1=#2\relax\dimen2=#3\relax\dimen5=#4%
\relax\vrule height\dimen5 depth\z@ width\z@\setbox8=\hbox to\dimen1{\tenln#1%
\hss}\cdna=\dimen5 \divide\cdna\dimen2 \ifnum\cdna=0 \then\box8 \else\dimen4=%
\dimen5 \advance\dimen4-\dimen2 \divide\dimen4\cdna\dimen3=\dimen1 \cdnb=%
\dimen2 \divide\cdnb1000 \divide\dimen3\cdnb\cdnb=\dimen4
\divide\cdnb1000 \multiply\dimen3\cdnb\dimen6\dimen1
\advance\dimen6-\dimen3 \cdnb=0
\ifPositiveGradient\then\dimen7\z@\else\dimen7\cdna\dimen4
\multiply\dimen4-1 \fi\loop\raise\dimen7\copy8
\ifnum\cdnb<\cdna\hskip-\dimen6 \advance\cdnb1
\advance\dimen7\dimen4
\repeat\fi}}\newdimen\objectheight\objectheight1.5ex

\def\execDiagonalLine{\setbox0=\hbox\tozpt{\cdna=\xcoord\cdnb=\ycoord\dimen8=%
\wd2 \dimen9=\ht2 \dimen0=\cdnb\DiagramCellHeight\advance\dimen0-2%
\MapShortFall\advance\dimen0-\objectheight\setbox2=\makeline\box2(\dimen8,%
\dimen9;.5\dimen0)\setbox4=\makeline\box4(\dimen8,\dimen9;.5\dimen0)\dimen0=2%
\wd2 \advance\dimen0-\cdna\DiagramCellWidth\advance\dimen0
2\DiagramCellWidth
\dimen2\DiagramCellHeight\advance\dimen2-\MapShortFall\dimen1\dimen2
\advance \dimen1-\ht1 \advance\dimen2-\ht2 \dimen6=\dimen2
\advance\dimen6.25\dimen8 \dimen3\dimen2 \advance\dimen3-\ht3
\dimen4=\dimen2 \dimen7=\dimen2 \advance \dimen4-\ht4
\advance\dimen7-\ht7 \advance\dimen7-.25\dimen8
\ifPositiveGradient\then\hss\raise\dimen4\hbox{\rlap{\box5}\box4}\llap{\raise
\dimen6\box6\kern.25\dimen9}\else\kern-.5\dimen0 \rlap{\raise\dimen1\box1}%
\raise\dimen2\box2 \llap{\raise\dimen7\box7\kern.25\dimen9}\fi\raise\dimen3%
\hbox\tozpt{\hss\box3\hss}\ifPositiveGradient\then\rlap{\kern.25\dimen9\raise
\dimen7\box7}\raise\dimen2\box2\llap{\raise\dimen1\box1}\kern-.5\dimen0
\else
\rlap{\kern.25\dimen9\raise\dimen6\box6}\raise\dimen4\hbox{\box4\llap{\box5}}%
\hss\fi}\ht0\z@\dp0\z@\box0}

\def\NorthWest{\PositiveGradientfalse\Climbingtrue\DiagonalChoice0 }\def
\NorthEast{\PositiveGradienttrue\Climbingtrue\DiagonalChoice1
}\def\SouthWest
{\PositiveGradienttrue\Climbingfalse\DiagonalChoice3 }\def\SouthEast{%
\PositiveGradientfalse\Climbingfalse\DiagonalChoice2 }

\newif\ifmoremapargs\def\gettwoargs{\setbox7=\box\voidb@x\setbox6=\box
\voidb@x\moremapargstrue\def\whichlabel{6}\def\xcoord{2}\def\ycoord{2}\def
\contgetarg{\def\whichlabel{7}\ifmoremapargs\then\let\next\getanarg\let
\contgetarg\execmap\else\let\next\execmap\fi\next}\getanarg}\def\getanarg{%
\futurelet\thenexttoken\switcharg}\def\getlabel#1#2#3{\setbox#1=\hbox{$%
\labelstyle\>{#3}\>$}\dimen0=\ht#1\advance\dimen0 .4ex\ht#1=\dimen0 \dimen0=%
\dp#1\advance\dimen0 .4ex\dp#1=\dimen0 \contgetarg}\def\eatspacerepeat{%
\afterassignment\getanarg\let\junk=
}\def\catcase#1:{{\ifcat\noexpand
\thenexttoken#1\then\global\let\xcase\docase\fi}\xcase}\def\tokcase#1:{{\ifx
\thenexttoken#1\then\global\let\xcase\docase\fi}\xcase}\def\default:{\docase}%
\def\docase#1\esac#2\esacs{#1}\def\skipcase#1\esac{}\def\getcoordsrepeat(#1,#%
2){\def\xcoord{#1}\def\ycoord{#2}\getanarg}\let\esacs\relax\def\switcharg{%
\global\let\xcase\skipcase\catcase{&}:\moremapargsfalse\contgetarg\esac
\catcase\bgroup:\getlabel\whichlabel-\esac\catcase^:\getlabel6\esac\catcase_:%
\getlabel7\esac\tokcase{~}:\getlabel3\esac\tokcase(:\getcoordsrepeat\esac
\catcase{
}:\eatspacerepeat\esac\default:\moremapargsfalse\contgetarg\esac
\esacs}

\catcode`\@=\cdna

\def\LH{{\bf LH}}\def\LM{{\bf LM}}\def\LT{{\bf
LT}}\def\KS{K\langle X\rangle}
\def\B{{\cal B}} \def\LC{{\bf LC}}
\def\G{{\cal G}} \def\FRAC#1#2{\displaystyle{\frac{#1}{#2}}}
\def\SUM^#1_#2{\displaystyle{\sum^{#1}_{#2}}} \def\O{{\cal O}}  \def\J{{\bf J}}
\def\RS{R\langle X\rangle} \def\SX{S\langle X\rangle}
\par

\section*{0. Introduction}
In the structure theory and the representation theory of associative
algebras over a ground field $K$, it is well known that numerous
popularly studied algebras have defining relations which form a
Gr\"obner basis in the classical sense (e.g., [Mor], [Gr]),  and
such algebras can be studied in a computational way via their
Gr\"obner defining relations (e.g., see [An], [CU], [GI-L], [Gr],
[Li2, 3], [Uf1, 2]); also we know that algebras defined by the
relations of the same type over a commutative ring $R$ are equally
important, for instance, the algebras over rings considered in
[Yam], [Ber], [CE], and [LVO2]. So, naturally we expect that certain
algebras over rings could be studied by means of Gr\"obner basis
theory as in loc. cit., and thus we hope that the following
statement would hold true:{\parindent=.6truecm\par

\item{$\bullet$ }  Let $K\langle X\rangle =K\langle
X_1,...,X_n\rangle$ be the free algebra of $n$ generators over a
field $K$, and let $R\langle X\rangle =R\langle X_1,...,X_n\rangle$
be the free algebra of  $n$ generators over an arbitrary commutative
ring $R$. By a Gr\"obner basis for an ideal in a free algebra we
mean the one as defined in [Mor] and [Gr]. If, with respect to some
monomial ordering $\prec$ on $\KS$, a subset $\G\subset\KS$ is a
Gr\"obner basis for the ideal $\langle\G\rangle$ in $\KS$, then,
taking a counterpart of $\G$ in $\RS$ (if it exists, again denoted
by $\G$) and using the same monomial ordering $\prec$ on $\RS$,
$\G$ is a Gr\"obner basis for the ideal $\langle\G\rangle$ in $\RS$.
\parindent=0pt\par

To see at what level the above statement may hold true, it is
necessary to see whether a version of the classical termination
theorem ([Mor], [Gr]) works well for verifying Gr\"obner bases in
$\RS$. }\par

Let $\KS =K\langle X_1,...,X_n\rangle$ be the free associative
$K$-algebra of $n$ generators over a field $K$, and let $\B$ be the
standard $K$-basis of $\KS$ consisting of monomials (words in
alphabet $X=\{ X_1,...,X_n\}$, including empty word which is
identified with the multiplicative identity element 1 of $\KS$).
Given a monomial ordering $\prec$ on $\B$ (i.e. a well-ordering
$\prec$ on $\B$ satisfying: $u\prec v$ implies $wus\prec wvs$ for
all $w,u,v,s\in\B$), and $f,g\in\KS-\{ 0\}$, if there are monomials
$u,v\in\B$ such that {\parindent=.75truecm\par
\item{(1)} $\LM (f)u=v\LM (g)$, and
\item{(2)} $\LM (f){\not |}~v$ and $\LM (g)\not |~u$,}{\parindent=0pt\par then the
element
$$o(f,u;~v,g)=\frac{1}{\LC (f)}(f\cdot u)-
\frac{1}{\LC (g)}(v\cdot g)$$ is referred to as an overlap element
of $f$ and $g$, where, with respect to $\prec$, $\LM (~)$ denotes
the function taking the leading monomial and $\LC (~)$ denotes the
function taking the leading coefficient on elements of $\KS$
respectively. Over the ground field $K$, the termination theorem in
the sense of ([Mor], [Gr]), which is known an algorithmic version of
Bergman's diamond lemma [Ber1], states
that}{\parindent=.5truecm\vskip 6pt
\item{$\bullet$} if ${\cal G}$ is an LM-reduced subset of $\KS$ (i.e.,
$\LM (g_i){\not |}~\LM (g_j)$ for $g_i,g_j\in\G$ with $i\ne j$),
then $\G$ is a Gr\"obner basis for the ideal $I=\langle\G\rangle$ if
and only if for each pair $g_i,g_j\in\G$, including $g_i=g_j$, every
overlap element $o(g_i,u;~v,g_j)$ of $g_i$ and $g_j$ has the
property $\overline{o(g_i,u;~v,g_j)}^{\G}=0,$ that is,
$o(g_i,u;~v,g_j)$ is reduced to 0 by division by $\G$;
\par}{\parindent=0pt\vskip 6pt and it follows that there is a
noncommutative analogue of the Buchberger algorithm for constructing
a (possibly infinite) Gr\"obner basis starting with a given finite
subset in $\KS$. Note that the algorithmic feasibility of the above
criterion lies in the fact that \parindent=.8truecm \par

\item{(a)} for each pair $(g_i,g_j)$ there are only finitely many
associated overlap elements, and

\item{(b)} there is no trouble with taking the inverse of a nonzero
coefficient when the division algorithm is performed, for, $K$ is a
field. \parindent=0pt\par

However, if the field $K$ is replaced by a commutative ring $R$,
and if $\G\subset\RS =R\langle X_1,...,X_n\rangle$ is taken such
that the leading coefficient $\LC (g)$ of some $g\in\G$ is not
invertible, then, even if $R$ is an arithmetic ring (e.g. the ring
$\mathbb{Z}$ of integers) as recently  considered by [Gol], there
seems no implementable  termination theorem (as we mentioned above)
for $\G$. Nevertheless, we have the following
observations:\parindent=.8truecm
\par

\item{(1)} If $\G$ is a Gr\"obner basis for an ideal $I$ in the free
algebra $\KS$ over a field $K$, then we may always assume that all
elements of $\G$ are monic, i.e., $\LC (g)=1$ for every $g\in\G$.
\item{(2)} If $S$ is a subset consisting of monic elements in the free algebra $\RS$
over a ring $R$, i.e., $\LC (f)=1$ for every $f\in S$, then, with
respect to any monomial ordering $\prec$ on $\RS$,  a division
algorithm by $S$ can be implemented in $\RS$ exactly as in
$\KS$.\parindent=0pt

Recalling the proof of the classical termination theorem ([Mor],
[Gr]), the above observations provide us with sufficient reason to
have an implementable termination theorem (as mentioned above) for
verifying whether certain given monic elements of  $R\langle
X\rangle$ form a Gr\"obner basis in the classical sense, so that our
foregoing expectation may come true. To present the details, we
organize this paper as follows. In Section 1, after giving a quick
introduction of the notion of a monic Gr\"obner basis in $\RS$, we
examine carefully that a version of the termination theorem in the
sense of ([Mor], [Gr]) holds true for verifying LM-reduced monic
Gr\"obner bases in $\RS$, just for convincing ourselves and also for
the reader's convenience from the viewpoint of ``to see is to
believe". This enables us to clarify that every monic Gr\"obner
basis in the free algebra $\KS$ over a field $K$ may give rise to a
monic Gr\"obner basis of the same type in the free algebra $\RS$,
and vice versa. In Section 2, by strengthening and generalizing a
result of [Li2], we show how PBW $R$-bases and monic Gr\"obner bases
of certain type can determine each other.  In the final Section 3,
we show that the working principle via PBW isomorphism developed in
[LWZ] and [Li3] can be generalized to study  quotient algebras of
$\RS$, so that many global structural properties of $R$-algebras
defined by monic Gr\"obner bases may be determined through their
$\NZ$-leading homogeneous algebras and $\B_R$-leading homogeneous
algebras.}\v5

Unless otherwise stated, rings considered in this paper are
associative rings with multiplicative identity 1,  ideals are meant
two-sided ideals, and modules are unitary left modules. For a subset
$U$ of a ring $S$, we write $\langle U\rangle$ (or ${~}_S\langle
U\rangle_S$ if necessary) for the ideal generated by $U$. Moreover,
we use $\NZ$, respectively $\mathbb{Z}$, to denote the set of
nonnegative integers, respectively the set of integers. \v5

\section*{1. Monic Gr\"obner Bases over $K$ vs Monic Gr\"obner Bases over $R$}
Let $R$ be an arbitrary commutative ring, $\RS =R\langle
X_1,...,X_n\rangle$ the free $R$-algebra of $n$ generators, and
$\B_R$ the standard $R$-basis of $\RS$ consisting of monomials
(words in alphabet $X=\{ X_1,...,X_n\}$, including empty word  which
is identified with the multiplicative identity element 1 of $\RS$).
Unless otherwise stated, monomials in $\B_R$ are denoted by lower
case letters $u,v,w,s,t,\cdots$.  In this section, after introducing
the notion of a monic Gr\"obner basis in $\RS$ and examining
carefully that a version of the termination theorem in the sense of
([Mor], [Gr]) holds true for verifying LM-reduced monic Gr\"obner
bases in $\RS$, we clarify that every monic Gr\"obner basis in the
free algebra $K\langle X_1,...,X_n\rangle$ over a field $K$ has a
counterpart in the free algebra $\RS$, and vice versa.\v5

First note that all monomial orderings used for free algebras over a
field can be well defined on the standard $R$-basis $\B_R$ of $\RS$.
In particular, by an $\NZ$-{\it graded monomial ordering} on $\B_R$,
denoted $\prec_{gr}$, we mean a monomial ordering on $\B_R$ which is
defined subject to a well-ordering $\prec$ on $\B_R$, that is,  for
$u,v\in\B_R$, $u\prec_{gr} v$ if either deg$u<$ deg$v$ or deg$u=$
deg$v$ but $u\prec v$, where deg$(~)$ denotes the degree function on
elements of $\RS$ with respect to a fixed {\it weight} $\NZ$-{\it
gradation} of $\RS$ (i.e. each $X_i$ is assigned  a positive degree
$n_i$, $1\le i\le n$). For instance, the usual $\NZ$-graded
(reverse) lexicographic ordering is a popularly used $\NZ$-graded
monomial ordering.\par

If $\prec$ is a monomial ordering on $\B_R$ and
$f=\sum_{i=1}^s\lambda_iw_i\in \RS$, where $\lambda_i\in R-\{ 0\}$
and $w_i\in\B_R$, such that $w_1\prec w_2\prec\cdots\prec w_s$, then
the {\it leading monomial} of $f$ is defined as $\LM (f)=w_s$ and
the {\it leading coefficient} of $f$ is defined as $\LC
(f)=\lambda_s$. For a subset $H\subset\RS$, we write $\LM (H)=\{\LM
(f)~|~f\in H\}$ for the set of leading monomials of $S$. We say that
a subset $G\subset\RS$ is {\it monic} if $\LC (g)=1$ for all $g\in
G$. Moreover, for $u,v\in\B_R$, as usual we say that $v$ divides
$u$, denoted $v|u$, if $u=wvs$ for some $w,~s\in\B_R$. \par With
notation and all definitions as above, it is easy to see that a
division algorithm by a monic subset $G$ is valid in $\RS$ with
respect to any fixed monomial ordering $\prec$ on $\B_R$. More
precisely, let $f\in\RS$. Noticing $\LC (g)=1$ for all $g\in G$, if
$\LM (g)|\LM (f)$ for some $g\in G$, then $f$ can be written as
$f=\LC (f)ugv+f_1$ with $u$, $v\in\B_R$, $f_1\in\RS$ satisfying $\LM
(f_1)\prec\LM (f)$; if $\LM (g){\not |}~\LM (f)$ for all $g\in G$,
then $f=f_1+\LC (f)\LM (f)$ with $f_1=f-\LC (f)\LM (f)$ satisfying
$\LM (f_1)\prec\LM (f)$. Next, consider the divisibility of $\LM
(f_1)$ by $\LM (g)$ with $g\in G$, and so forth. Since $\prec$ is a
well-ordering, after a finite number of successive division by
elements in $G$ in this way, we see that $f$ can be written as
$$\begin{array}{rcl} f&=&\sum_{i,j}\lambda_{ij}u_{ij}g_jv_{ij}+r_f,~\hbox{where}~\lambda_{ij}\in R,~u_{ij},v_{ij}\in\B_R,~g_j\in G,\\
&{~}&\hbox{and}~r_f=\sum_p\lambda_pw_p~\hbox{with}~\lambda_p\in
R,~w_p\in\B_R,\\
&{~}&\hbox{satisfying}~\LM (u_{ij}g_jv_{ij})\preceq\LM
(f)~\hbox{whenever}~\lambda_{ij}\ne 0,\\
&{~}&\LM (r_f)\preceq\LM (f)~\hbox{and}~\LM (g)\not |~w_p~\hbox{for
every}~g\in G~\hbox{whenever}~\lambda_p\ne 0.\end{array}$$
If, in the presentation of $f$ above, $r_f=0$, then we say that $f$
{\it is reduced to} 0 {\it by division by} $G$, and we write $\OV
f^G=0$ for this property.\par
The validity of such a division algorithm by $G$ leads to the
following definition.{\parindent=0pt\v5

{\bf 1.1. Definition} Let $\prec$ be a fixed monomial ordering on
$\B_R$, and $I$ an ideal of $\RS$. A {\it monic Gr\"obner basis} of
$I$ is a subset $\G\subset I$ satisfying:\par (1) $\G$ is monic; and
\par                                     (2) $f\in I$ and $f\ne 0$
implies $\LM (g)|\LM (f)$ for some $g\in \G$. }\v5 By the division
algorithm presented above, it is clear that a monic Gr\"obner basis
of $I$ is first of all a generating set of the ideal $I$, and
moreover, a monic Gr\"obner basis of $I$ can be characterized as
follows. {\parindent=0pt\v5

{\bf 1.2. Proposition}  Let $\prec$ be a fixed monomial ordering on
$\B_R$, and $I$ an ideal of $\RS$. For a monic subset $\G\subset I$,
the following statements are equivalent:\par (i) $\G$ is a monic
Gr\"obner basis of $I$;\par (ii) Each nonzero $f\in I$ has a
Gr\"obner representation:
$$\begin{array}{rcl} f&=&\sum_{i,j}\lambda_{ij}u_{ij}g_jv_{ij},
~\hbox{where}~\lambda_{ij}\in R,~u_{ij},v_{ij}\in\B_R,~g_j\in G,\\
&{~}&\hbox{satisfying}~\LM (u_{ij}g_jv_{ij})\preceq\LM
(f)~\hbox{whenever}~\lambda_{ij}\ne 0,\end{array}$$ or equivalently,
$\OV f^{\G}=0$; \par (iii) $\langle\LM (\G )\rangle =\langle\LM
(I)\rangle$. \QED}\v5

Let $\prec$ be a monomial ordering on the standard $R$-basis $\B_R$
of $\RS$, and let $G$ be a monic subset of $\RS$. We call an element
$f\in \RS$ a {\it normal element} (mod $G$) if $f=\sum_j\mu_jv_j$
with $\mu_j\in R$, $v_j\in\B_R$, and $f$ has the property that $\LM
(g)\not |~v_j$ for every $g\in G$ and all $\mu_j\ne 0$. The set of
normal monomials in $\B_R$ (mod $G$) is denoted by $N(G)$, i.e.,
$$N(G)=\{ u\in\B_R~|~\LM (g){\not |}~u, ~g\in G\} .$$
Thus, an element $f\in\RS$ is normal (mod $G$) if and only if $f\in
\sum_{u\in N(G)}Ru$. {\parindent=0pt\v5

{\bf 1.3. Proposition} Let $\G$ be a monic Gr\"obner basis of the
ideal $I=\langle\G\rangle$ in $\RS$ with respect to some monomial
ordering $\prec$ on $\B_R$. Then each nonzero $f\in \RS$ has a
finite presentation
$$f=\sum_{i,j}\lambda_{ij}s_{ij}g_iw_{ij}+r_f,\quad \lambda_{ij}\in R,~s_{ij},
w_{ij}\in\B_R,~g_i\in \G,$$ where $\LM (s_{ij}g_iw_{ij})\preceq\LM
(f)$ whenever $\lambda_{ij}\ne 0$, and either $r_f=0$ or $r_f$ is a
unique normal element (mod $\G$). Hence, $f\in I$ if and only if
$r_f=0$, solving the ``membership problem" for $I$.\par \vskip 6pt
{\bf Proof} By the division algorithm by $\G$, $f$ can be written as
$f=\sum_{i,j}\lambda_{ij}s_{ij}g_iw_{ij}+r_f$ where either $r_f=0$
or $r_f$ is normal. Suppose after division by $\G$ we also have
$f=\sum_{t,j}\lambda_{tj}s_{tj}g_tw_{tj}+r$, where $r$ is normal
(mod $\G$). Then $r-r_f\in I$ and hence there is some $g\in\G$ such
that  $\LM (g)|\LM (r-r_f)$. But by the definition of a normal
element this is possible only if $r=r_f$.\QED} \v5

The foregoing discussion enables us to obtain further
characterization of a monic Gr\"obner basis $\G$, which, in turn,
gives rise to the  fundamental decomposition theorem of $\RS$
(viewed as an $R$-module) by the ideal $I=\langle\G\rangle$, and
gives rise to a free $R$-basis for the $R$-algebra $\RS /I$.
{\parindent=0pt\v5

{\bf 1.4. Theorem}  Let  $I=\langle\G\rangle$ be an ideal of $\RS$
generated by a monic subset $\G$. With notation as above, the
following statements are equivalent.\par (i) $\G$ is a monic
Gr\"obner basis of $I$. \par (ii) The $R$-module $\RS$ has the
decomposition
$$\RS=I\oplus \sum_{u\in N(\G )}Ru=\langle\LM (I)\rangle\oplus\sum_{u\in N(\G )}Ru.$$
(iii) The canonical image $\OV{N(\G )}$ of $N(\G )$ in $\RS
/\langle\LM (I)\rangle$ and $\RS/I$ forms a free $R$-basis for $\RS
/\langle\LM (I)\rangle$ and $\RS/I$ respectively. \par\QED}\v5

With notation and every definition introduced so far, we proceed now
to show that a version of the termination theorem in the sense of
([Mor], [Gr]) holds true for verifying an LM-reduced monic Gr\"obner
basis in $\RS$ (see the definition below). \v5

Given a monomial ordering $\prec$ on $\B_R$, we say that a subset
$G\subset\RS$ is {\it LM-reduced} if
$$\LM (g_i)\not |~\LM (g_j)~\hbox{for all}~g_i,g_j\in G~\hbox{with}~g_i\ne g_j.$$
If a subset $G\subset\RS$ is both LM-reduced and monic, then we call
$G$ an {\it LM-reduced monic subset}. Thus we have the notion of an
{\it LM}-{\it reduced monic Gr\"obner basis}.\par
Let $I$ be an ideal of $\RS$. If $\G$ is a monic Gr\"obner basis of
$I$ and $g_1,g_2\in\G$ such that $g_1\ne g_2$ but $\LM (g_1)|\LM
(g_2)$, then clearly $g_2$ can be removed from $\G$ and the remained
subset $\G-\{ g_2\}$ is again a monic Gr\"obner basis for $I$.
Hence, in order to have a better criterion for monic Gr\"obner basis
we need only to consider the subset which is both LM-reduced and
monic. \v5

Let $\prec$ be a monomial ordering on $\B_R$. For two monic elements
$f,~g\in\RS-\{ 0\}$, including $f=g$, if there are monomials
$u,v\in\B_R$ such that {\parindent=.75truecm\par
\item{(1)} $\LM (f)u=v\LM (g)$, and
\item{(2)} $\LM (f)\not |~v$ and $\LM (g)\not |~u$,}{\parindent=0pt\par then the
element
$$o(f,u;~v,g)=f\cdot u-
v\cdot g$$ is called an {\it overlap element} of $f$ and $g$. From
the definition we are clear about the fact that
$$\LM ((o(f,u;~v,g))\prec \LM (fu)=\LM (vg),$$
and moreover, there are only finitely many overlap elements for each
pair $(f,g)$ of elements in $\RS$.}\par

With the preparation made above, below we  mention a termination
theorem for checking LM-reduced monic Gr\"obner bases in $\RS$, and,
also we present a careful step-by-step verification of its
correctness for the reason that we are working on an arbitrary ring
instead of a field after all, though the process is only a light
modification of the argument given in [Gr].  {\parindent=0pt\v5

{\bf 1.5. Theorem} (Termination theorem)  Let $\prec$ be a fixed
monomial ordering on $\B_R$. If $\G$ is an LM-reduced monic subset
of $\RS$, then $\G$ is an LM-reduced monic Gr\"obner basis for the
ideal $I=\langle\G\rangle$ if and only if for each pair
$g_i,g_j\in\G$, including $g_i=g_j$, every overlap element
$o(g_i,u;~v,g_j)$ of $g_i$, $g_j$ has the property
$\overline{o(g_i,u;~v,g_j)}^{\G}=0,$ that is, by division by $\G$,
every $o(g_i,u;~v,g_j)$ is reduced to zero.\vskip 6pt

{\bf Proof} Since  $o(g_i,u;~v,g_j)\in I$, the necessity follows
from Proposition 1.2.}\par Under the assumption on overlap elements
we prove the sufficiency by showing that if $f\in I$ then $\LM
(g)|\LM (f)$ for some $g\in\G$. Suppose the contrary that $\LM
(g){\not |}~\LM (f)$ for any $g\in\G$. Then we proceed to derive a
contradiction.
\par Since $I=\langle\G\rangle$, $f$ may be presented
as a finite sum
$$f=\sum_{i,j}\lambda_{ij}v_{ij}g_iw_{ij},~
\hbox{where}~\lambda_{ij}\in R,~v_{ij},w_{ij}\in\B_R,~\hbox{and}~
g_i\in\G.\leqno{(1)}$$ Let $u$ be the largest monomial occurring on
the right hand side of (1). Then noticing that each $g_i$ is monic,
$u$ occurs as some $v_{ij}\LM (g_i)w_{ij}$. Since $\LM (g_i){\not
|}~\LM (f)$ for the $g_i$ occurring in (1), it follows that $\LM
(f)\prec u$ and $u$ must occur at least twice on the right hand side
of (1) for a cancelation, that is, we may have
$$u=v_{ij}\LM (g_i)w_{ij}=v_{k\ell}\LM (g_k)w_{k\ell}.\leqno{(2)}$$
Among all such presentations of $f$ we can choose one such that
\vskip 6pt{\parindent=.7truecm
\item{(3)} $u$ has the fewest occurrences
on the right hand side of (1) and $u$ is as small as possible.\vskip
6pt\parindent=0pt To go further, let us simplify notation by writing
$v=v_{ij}$, $g=g_i$, $w=w_{ij}$, $v'=v_{k\ell}$, $g'=g_k$, and
$w'=w_{k\ell}$. Thus, the above (2) is turned into the form
$$u=v\LM (g)w=v'\LM (g')w'.\leqno{(4)}$$
Moreover, as usual we use $l(s)$ to denote the length of a monomial
$s\in\B_R$. Below we show, through a case by case study of the above
(4), that the choice of $f$ satisfying the above (3) is
impossible.\par {\parindent=0pt\par {\bf Case A}:
$l(v)<l(v')$.}{\parindent=0pt\par {\bf Case A.1}: $l(w)<l(w')$.}\par
This implies that $\LM (g)$ contains $\LM (g')$ as a subword, and
hence, $\LM (g')|\LM (g)$, contradicting the hypothesis on
$\G$.{\parindent=0pt\par {\bf Case A.2}: $l(w)\ge l(w')$. Then we
have to deal with two possibilities.\par {\bf Case A.2.1}: $l(v')\ge
l(v\LM (g))$.}\par This implies that there is no overlap element of
$\LM (g)$ and $\LM (g')$ in $u$. By the assumption on lengths, it
follows that there is a segment $w''$ of $v'$ such that $v'=v\LM
(g)w''$ and $w=w''\LM (g')w'$, i.e.,
$$u=v\LM (g)w''\LM (g')w'.$$
Rewriting $g$ as $g=\LM (g)+\sum\lambda_iu_i$, $g'=\LM
(g')+\sum\mu_iu_i'$, then
$$\begin{array}{rcl} vgw&=&vgw''g'w'-vgw''
\left (g'-\LM (g')\right )w'\\
&=&v\LM (g)w''g'w'+{\sum}\lambda_ivu_iw''g'w'
-{\sum}\mu_ivgw''u_i'w'\\
&=&v'g'w'+{\sum}\lambda_ivu_iw''g'w'
-{\sum}\mu_ivgw''u_i'w'.\end{array}$$ Thus, in writing $vgw$ this
way, the number of occurrences of $u$ in the chosen presentation of
$f$ satisfying the above (3) can be reduced, a
contradiction.{\parindent=0pt\par {\bf Case A.2.2}: $l(v')<l(v\LM
(g))$.}\par This implies that there is an overlap element of $\LM
(g)$ and $\LM (g')$ in $u$, that is, there is a segment $r$ of $w$
and a segment $s$ of $v'$ such that $vs=v'$, $rw'=w$ and $\LM
(g)r=s\LM (g')$. Hence, $$u=v\LM (g)rw'=vs\LM (g')w'~\hbox{and}~
o(g,r;~s,g')=gr- sg'.\leqno{(5)}$$ Furthermore, it follows from
$gr=o(g,r;~s,g')+sg'$ that
$$vgw=vgrw'=v\cdot o(g,r;~s,g')\cdot w'+v'g'w'.\leqno{(6)}$$
By the assumption, $o(g,r;~s,g')$ is reduced to 0 under the division
by $\G$, namely,
$$o(g,r;~s,g')=\sum_{k,j}c_{kj}v_{kj}g_kw_{kj},~v_{kj},w_{kj}\in\B_R,
~g_k\in\G,~c_{kj}\in R,\leqno{(7)}$$ satisfying
$$\hbox{if}~c_{kj}\ne 0~\hbox{then}~v_{kj}\LM (g_k)w_{kj}\prec \LM (g\cdot r)=\LM (s\cdot g').\leqno{(8)}$$
Combining the above (5) -- (8), once again we see that the number of
occurrences of $u$ in the chosen presentation of $f$ satisfying the
foregoing (3) can be reduced, a contradiction. {\parindent=0pt\par
{\bf Case B}: $l(v)=l (v')$.}\par This implies $\LM (g)|\LM (g')$ or
$\LM (g')|\LM (g)$, which contradicts the assumption that $\G$ is
LM-reduced.{\parindent=0pt\par {\bf Case C}: $l(v)>l(v')$.}\par This
is similar to Case 1.\QED\v5

{\bf Remark} (i) Let $\G$ be an LM-reduced monic subset of $\RS$ and
$I=\langle\G\rangle$. Since for each pair $g_i,~g_j\in\G$, including
$g_i=g_j$, there are only finitely many associated overlap elements,
by Theorem 1.5 we can check effectively whether $\G$ is a Gr\"obner
basis of $I$ or not, when $\G$ is a finite subset.\par

(ii) Obviously, if $\G\subset\RS$ is an LM-reduced subset with the
property that each $g\in \G$ has the leading coefficient $\LC (g)$
which is invertible in $R$, then Theorem 1.5 is also valid for
$\G$.\par

(iii) It is obvious as well that Theorem 1.5 does not necessarily
induce an analogue of the Buchberger algorithm as in the classical
case. \par

(iv) It is not difficult to see that all discussion we presented so
far is valid for getting monic Gr\"obner bases in a commutative
polynomial ring $R[x_1,...,x_n]$ over an arbitrary commutative ring
$R$ where overlap elements are replaced by $S$-polynomials.} \v5

Note that, throughout the proof of Theorem 1.5, nothing involves the
invertibility of an element in the ring $R$; moreover, division
algorithm by monic elements in a free algebra (over a field or over
a commutative ring) never touches on the invertibility of a
coefficient. Therefore, we are now clear about the relation between
monic Gr\"obner bases over a field and monic Gr\"obner bases over a
commutative ring.{\parindent=0pt\v5

{\bf 1.6. Proposition} Let $K\langle X\rangle =K\langle
X_1,...,X_n\rangle$ be the free algebra of $n$ generators over a
field $K$, and let $R\langle X\rangle =R\langle X_1,...,X_n\rangle$
be the free algebra of  $n$ generators over an arbitrary commutative
ring $R$. With notation as before, fixing the same monomial ordering
$\prec$ on both $\KS$ and $\RS$, the following statements hold.\par
(i) If a monic subset $\G\subset\KS$ is a Gr\"obner basis for the
ideal $\langle\G\rangle$ in $\KS$, then, taking a counterpart of
$\G$ in $\RS$ (if it exists), again denoted by $\G$, $\G$ is a monic
Gr\"obner basis for the ideal $\langle\G\rangle$ in $\RS$. \par

(ii) If a monic subset $\G\subset\RS$ is a Gr\"obner basis for the
ideal $\langle\G\rangle$ in $\RS$, then, taking a counterpart of
$\G$ in $\KS$ (if it exists), again denoted by $\G$, $\G$ is a
Gr\"obner basis for the ideal $\langle\G\rangle$ in $\KS$.\par
\QED}\v5

From the literature we know that numerous well-known $K$-algebras
over a field $K$, such as the $n$-th Weyl algebra over $K$, the
enveloping algebra of a $K$-Lie algebra, a $K$-exterior algebra, a
$K$-Clifford algebra, a down-up $K$-algebra, etc., all have defining
relations that form an LM-reduced monic  Gr\"obner basis in a free
$K$-algebra. Hence, by Proposition 1.6, if the field $K$ is replaced
by a commutative ring $R$, then all of these $R$-algebras have
defining relations that form an LM-reduced monic Gr\"obner basis in
a free $R$-algebra. Below we give another example illustrating
Theorem 1.5 and Proposition 1.6.{\parindent=0pt\v5

{\bf Example 1.} Let $R$ be a commutative ring. Consider in $\RS
=R\langle X_1,...,X_n\rangle$ the subset $\G =\Omega\cup {\cal R}$
consisting of
$$\begin{array}{l} \Omega\subseteq\{ g_i=X_i^p~|~1\le i\le n\}~\hbox{with}~
p\ge 2~\hbox{a fixed integer},\\
{\cal R}=\{ g_{ji}=X_jX_i-\lambda_{ji}X_iX_j~|~1\le i<j\le
n\}~\hbox{with}~\lambda_{ji}\in R,\\
\quad\quad~\hbox{that is,}~\lambda_{ji}~\hbox{may be
zero}.\end{array}$$  In the case that $R=K$ is a field, it was
verified in ([Li4], Example 4) that, under the $\NZ$-graded
lexicographic ordering $\prec_{gr}$ such that
$X_1\prec_{gr}X_2\prec_{gr}\cdots\prec_{gr} X_n$, $\G$ forms an 
LM-reduced monic Gr\"obner basis for the ideal $I=\langle \G\rangle$ 
in $\KS$. Hence, by Proposition 1.6, $\G$ is an LM-reduced monic 
Gr\"obner basis for the ideal $I=\langle\G \rangle$ in $\RS$. 
Furthermore, the division by $\LM (\G )$ yields
$$ N(\G )=\left\{ X_1^{\alpha_1}X_2^{\alpha_2}\cdots X_n^{\alpha_n}~\Big |~\alpha_i\in\NZ~\hbox{and}~0\le \alpha_s\le p-1~\hbox{if}~X^p_s\in\Omega\right\} .$$
It follows from Theorem 1.4 that both the algebras $\RS /I$ and $\RS
/\langle\LM (I)\rangle$ have the free $R$-basis
$$ \OV{N(\G )}=\left\{ \OV{X}_1^{\alpha_1}\OV{X}_2^{\alpha_2}\cdots \OV{X}_n^{\alpha_n}~\Big
|~\alpha_i\in\NZ~\hbox{and}~0\le \alpha_s\le
p-1~\hbox{if}~X^p_s\in\Omega\right\} ,$$ where each $\OV{X}_i$ is
the canonical image of $X_i$ in $\RS /I$ and $\RS /\langle\LM
(I)\rangle$ respectively.}\par

Here let us point out that this example contains two families of
special $R$-algebras, that is, in the case where $\Omega
=\emptyset$, the $R$-algebra $\RS /I$ is similar to the coordinate
ring of a quantum affine $n$-space over a field (such a quantum
coordinate ring over a field is defined with all the
$\lambda_{ji}\ne 0$); and in the case where $\Omega =\{
g_i=X_i^2~|~1\le i\le n\}$, the algebra $\RS /I$ is  similar to the
quantum grassmannian (or quantum exterior) algebra over a field in
the sense of [Man] (such a quantum grassmannian algebra over a field
is defined with all the $\lambda_{ji}\ne 0$).\v5

We finish this section by a useful corollary of Theorem
1.5.{\parindent=0pt\v5

{\bf 1.7. Corollary} Let $R$ be a commutative ring and $R'$ a
subring of $R$ with the same identity element 1. If we consider the
free $R$-algebra $\RS =R\langle X_1,...,X_n\rangle$ and the free
$R'$-algebra $R'\langle X\rangle =R'\langle X_1,...,X_n\rangle$,
then the following two statements are equivalent for a subset
$\G\subset R'\langle X\rangle$:\par

(i) $\G$ is an LM-reduced monic Gr\"obner basis for the ideal
$I=\langle\G\rangle$ in $R'\langle X\rangle$ with respect to some
monomial ordering $\prec$ on the standard $R'$-basis $\B_{R'}$ of
$R'\langle X\rangle$;\par

(ii) $\G$ is an LM-reduced monic Gr\"obner basis for the ideal
$J=\langle\G\rangle$ in $\RS$ with respect to the monomial ordering
$\prec$ on the standard $R$-basis $\B_R$ of $\RS$, where $\prec$ is
the same monomial ordering used in (i).\vskip 6pt {\bf Proof} Let
$\G\subset R'\langle X\rangle$ be an LM-reduced monic subset.
Noticing that $B_R=B_{R'}$, each pair $g_i,g_j\in\G$ has the same
set of overlap elements in both $R'\langle X\rangle$ and $\RS$.
Moreover, noticing that performing the division of an overlap
element $o(g_i,u;~v,g_j)$ by $\G$ in both $R'\langle X\rangle$ and
$\RS$ uses only coefficients from $R'$. It follows that
$\OV{o(g_i,u;~v,g_j)}^{\G}=0$ in $R'\langle X\rangle$ if and only if
$\OV{o(g_i,u;~v,g_j)}^{\G}=0$ in $\RS$. Hence the equivalence of (i)
and (ii) is proved by the termination theorem for LM-reduced monic
Gr\"obner bases in $R'\langle X\rangle$ and the termination theorem
for LM-reduced monic Gr\"obner bases in $\RS$, respectively.\QED
}\v5

\section*{2. PBW $R$-bases vs Specific Monic Gr\"obner Bases}
Let $R$ be a commutative ring and $A=R[a_1,...,a_n]$ a finitely
generated $R$-algebra with generators $a_1,...,a_n$. If the set
$\mathscr{B}=\{ a_1^{\alpha_1}a_2^{\alpha_2}\cdots
a_n^{\alpha_n}~|~\alpha_j\in\NZ\}$ forms a free $R$-basis of $A$,
that is, $A$ is, as an $R$-module, free with the basis
$\mathscr{B}$, then, in honor of the classical PBW
(Poincar\'e-Birkhoff-Witt) theorem for enveloping algebras of Lie
algebras over a ground field $K$, the set $\mathscr{B}$ is usually
referred to as a {\it PBW $R$-basis} of $A$. Presenting $A$ as a
quotient algebra of the free $R$-algebra $\RS =R\langle
X_1,...,X_n\rangle$, i.e., $A=\RS /I$ with $I$ an ideal of $\RS$,
the aim of this section is to show, under a mild condition, that $A$
has a PBW $R$-basis is equivalent to that $I$ has a specific monic
Gr\"obner basis.  This result enables us to obtain PBW $R$-bases by
means of monic Gr\"obner bases on one hand; and on the other hand,
since it is well known that in practice there are different ways to
find a PBW basis of a given algebra provided it exists (e.g., [Ros],
[Yam], [Rin], [Ber]), this result also enables us to obtain monic
Gr\"obner bases via already known PBW $R$-bases. \par

Throughout this section, we let $R$ be a commutative ring, $\RS
=R\langle X_1,...,X_n\rangle$ the free $R$-algebra of $n$
generators, and $\B_R$ the standard $R$-basis of $\RS$. All
notations and notions concerning monic Gr\"obner bases in $\RS$ are
maintained as before. \v5

Let $I$ be an ideal of $\RS$ such that the $R$-algebra $A=\RS /I$
has the PBW $R$-basis $\mathscr{B} =\left\{
\OV{X}_1^{\alpha_1}\OV{X}_2^{\alpha_2}\cdots
\OV{X}_n^{\alpha_n}~|~\alpha_i\in\NZ\right\}$, where each $\OV X_i$
is the canonical image of $X_i$ in $A$. Then $I$ contains
necessarily a subset $G$ consisting of $\frac{n(n-1)}{2}$ elements
of the form:
$$g_{ji}=X_jX_i-\sum_{\alpha}\lambda_{\alpha}w_{\alpha},~\hbox{where}~1\le i<j\le n,~
\lambda_{\alpha}\in
R,~w_{\alpha}={X}_1^{\alpha_1}{X}_2^{\alpha_2}\cdots
{X}_n^{\alpha_n}.$$

In light of Theorem 1.4 and the observation made above,  below we
give the main result of this section which, indeed, strengthens and
generalizes ([Li2], CH.III, Theorem 1.5). {\parindent=0pt\v5

{\bf 2.1. Theorem} Let $I$ be an ideal of $\RS$, $A=\RS /I$. Suppose
that $I$ contains a monic subset of $\frac{n(n-1)}{2}$ elements
$G=\{ g_{ji}~|~1\le i<j\le n\}$ such that, with respect to some
monomial ordering $\prec$ on the standard $R$-basis $\B_R$ of $\RS$,
$\LM (g_{ji})=X_jX_i$ for  $1\le i<j\le n$. The following two
statements are equivalent.
\par
(i) The $R$-algebra $A$ has the  PBW $R$-basis $\mathscr{B}=\{ \OV
X_1^{\alpha_1}\OV X_2^{\alpha_2}\cdots \OV
X_n^{\alpha_n}~|~\alpha_j\in\NZ\}$ where each $\OV X_i$ is the
canonical image of $X_i$ in $A$.\par (ii) Any monic subset $\G$ of
$I$ containing $G$ is a monic  Gr\"obner basis for $I$ with respect
to $\prec$. \vskip 6pt
{\bf Proof} (i) $\Rightarrow$ (ii) Let $\G$ be a monic  subset of
$I$ containing $G$, and let $$N(\G )=\{ u\in\B_R~|~\LM (g)\not
|~u,~g\in\G\}$$ be the set of normal monomials in $\B_R$ (mod $\G$).
If $f\in I$ and $f\ne 0$, then, after implementing the division of
$f$ by $\G$ (with respect to the given monomial ordering $\prec$) we
have
$$\begin{array}{rcl} f&=&\sum_{i,j}\lambda_{ij}u_{ij}g_iv_{ij}+r_f,~
\hbox{where}~\lambda_{ij}\in R,~u_{ij},v_{ij}\in\B_R,~g_i\in\G ,\\
&{~}&\hbox{satisfying}~\LM (w_{ij}g_iv_{ij})\preceq\LM
(f)~\hbox{whenever}~\lambda_{ij}\ne 0,\\
&{~}&\hbox{and}~r_f=\sum_p\lambda_pw_p~\hbox{with}~\lambda_p\in
R~\hbox{and}~w_p\in N(\G ).\end{array}$$ Note that $g_{ji}\in
G\subseteq\G$ and $\LM (g_{ji})=X_jX_i$ by the assumption.  It
follows that $N(\G )\subseteq\{ X_1^{\alpha_1}X_2^{\alpha_2}\cdots
X_n^{\alpha_n}~|~\alpha_j\in\NZ\}$. Thus, since $\mathscr{B}$ is a
free  $R$-basis of $A$,
$r_f=\sum_p\lambda_pw_p=f-\sum_{i,j}u_{ij}g_iv_{ij}\in I$ implies
$\lambda_p=0$ for all $p$. Consequently $r_f=0$. This shows that
every nonzero element of $I$ has a Gr\"obner representation by the
elements of $\G$. Hence $\G$ is a monic Gr\"obner basis for $I$ by
Proposition 1.2.\par (ii) $\Rightarrow$ (i) By (ii), the subset $G$
itself is a monic  Gr\"obner basis of $I$ with respect to $\prec$.
Let $N(G)$ be the set of normal monomials in $\B_R$ (mod $G$).
Noticing that $\LM (g_{ji})=X_jX_i$ for every $g_{ji}\in G$, it
follows that $N(G)= \{ X_1^{\alpha_1}X_2^{\alpha_2}\cdots
X_n^{\alpha_n}~|~\alpha_j\in\NZ\}$, and thereby the algebra $A$ has
the desired PBW $R$-basis $\mathscr{B}$ by Theorem 1.4.\QED}\v5

We illustrate Theorem 2.1 by several examples. The first four
examples given below serve to obtain  monic Gr\"obner bases by means
of already known PBW $R$-bases which are obtained in the literature
without using the theory of Gr\"obner basis.{\parindent=0pt\v5

{\bf Example 1.} (This is a special case of Example 3 given later.)
Let $\textsf{g}=R[V]$ be the $R$-Lie algebra defined by the free
$R$-module $V=\oplus_{i=1}^nRx_i$ and the bracket product
$[x_j,x_i]=\sum_{\ell =1}^n\lambda_{ji}^{\ell}x_{\ell}$, $1\le
i<j\le n$, $\lambda_{ji}^{\ell}\in R$. By the classical PBW theorem,
the universal enveloping algebra $U(\textsf{g})$ of $\textsf{g}$ has
the PBW $R$-basis $\mathscr{B}=\{ x_1^{\alpha_1}x_2^{\alpha_2}\cdots
x_n^{\alpha_n}~|~\alpha_j\in\NZ\}$. If, with respect to the natural
$\NZ$-gradation of $\RS =R\langle X_1,...,X_n\rangle$, we use an
$\NZ$-graded monomial ordering $\prec_{gr}$ on the standard
$R$-basis $\B_R$ of $\RS$  such that
$X_1\prec_{gr}X_2\prec_{gr}\cdots\prec_{gr} X_n$  (i.e., deg$X_i=1$,
$1\le i\le n$), then the set of defining relations
$$\G =\left\{\left. g_{ji}=X_jX_i-X_iX_j-\sum_{\ell =1}^n\lambda_{ji}^{\ell}X_{\ell}~\right |~
1\le i<j\le n\right\}$$ of  $U(\textsf{g})$ satisfies $\LM
(g_{ji})=X_jX_i$ for $1\le i<j\le n$. Hence, by Theorem 2.1, $\G$ is
a monic Gr\"obner basis for the ideal $I=\langle\G\rangle$ in $\RS$.
\v5

{\bf Example 2.} Let $U_q^+(A_N)$ be the $(+)$-part of the
Drinfeld-Jimbo quantum group of type $A_N$ over a commutative ring
$R$, where $q\in R$ is invertible and $q^8\ne 1$. This example shows
that the defining relations (Jimbo relations) of $U_q^+(A_N)$ over
$R$ form a monic Gr\"obner basis in a free $R$-algebra. By
Proposition 1.6, we reach this property over a field $K$. }\par

In [Ros] and [Yam] it was proved that, over a field $K$,
$U^+_q(A_N)$ has a PBW $K$-basis with respect to the defining
relations (Jimbo relations) of $U^+_q(A_N)$;  later in [BM] such a
PBW basis was recaptured by showing that the Jimbo relations form a
Gr\"obner basis ([BM] Theorem 4.1), where the proof was sketched to
check that all compositions (overlaps) of Jimbo relations reduce to
zero on the base argument of [Yam]. Recently, a very detailed
elementary verification of the fact that all compositions (overlaps)
of Jimbo relations reduce to zero and hence the Jimbo relations form
a Gr\"obner basis (namely Theorem 4.1 of [BM]) was carried out by
[CSS]. Now,  by using Theorem 2.1 we will see that it is indeed very
easy to conclude:  the Jimbo relations form a Gr\"obner basis. \par

Recall that the Jimbo relations (as described in [Yam]) are given by
$$\begin{array}{ll} x_{mn}x_{ij}-q^{-2}x_{ij}x_{mn},&((i,j),(m,n))\in C_1\cup C_3,\\
x_{mn}x_{ij}-x_{ij}x_{mn},&((i,j),(m,n))\in C_2\cup C_6,\\
x_{mn}x_{ij}-x_{ij}x_{mn}+(q^2-q^{-2})x_{in}x_{mj},&((i,j),(m,n))\in
C_4,\\
x_{mn}x_{ij}-q^2x_{ij}x_{mn}+qx_{in},&((i,j),(m,n))\in C_1\cup
C_3,\end{array}$$ where with $\Lambda_N=\{ (i,j)\in\NZ\times
\NZ~|~1\le i<j\le N+1\}$,
$$\begin{array}{ll} C_1=\{ ((i,j),(m,n))~|~i=m<j<n\},&C_2=\{ ((i,j),(m,n))~|~i<m<n<j\},\\
C_3=\{ ((i,j),(m,n))~|~i<m<j=n\},&C_4=\{ ((i,j),(m,n))~|~i<m<j<n\},\\
C_5=\{ ((i,j),(m,n))~|~i<j=m<n\},&C_6=\{
((i,j),(m,n))~|~i<j<m<n\}.\end{array}$$ By [Yam], for $q^8\ne 1$,
$U^+_q(A_N)$ has the PBW basis consisting of elements
$$x_{i_1j_1}x_{i_2j_2}\cdots x_{i_kj_k}~\hbox{with}~(i_1,j_1)\le (i_2,j_2)\le\cdots\le
(i_k,j_k),~k\ge 0,$$ where $(i_{\ell},j_{\ell})\in\Lambda_N$ and $<$
is the lexicographic ordering on $\Lambda_N$. If we use the
$\NZ$-graded monomial ordering $\prec_{gr}$ (on the standard
$K$-basis $\B$ of the corresponding free algebra) subject to
$$x_{ij}\prec_{gr} x_{mn}\Longleftrightarrow (i,j)<(m,n),$$ then it is
clear that for each pair $((i,j),(m,n))\in C_i$ with $(i,j)<(m,n)$,
the leading monomial of the corresponding relation is of the form
$x_{mn}x_{ij}$ as required by Theorem 2.1. \v5 {\parindent=0pt

{\bf Example 3.} With $\RS =R\langle X_1,...,X_n\rangle$, where $R$
is an arbitrary commutative ring, recall from [Ber] that a
$q$-algebra $A=\RS /\langle\G\rangle$ over $R$ is defined by the set
$\G$ of quadric relations
$$\begin{array}{lll}
g_{ji}&=&X_jX_i-q_{ji}X_iX_j-\{ X_j,X_i\},~1\le i<j\le
n,~\hbox{where}~q_{ji}\in R-\{0\},\\
&{~}&\hbox{and}~\{ X_j,X_i\}
=\displaystyle{\sum}\alpha^{k\ell}_{ji}X_kX_{\ell}+
\displaystyle{\sum}\alpha_hX_h+c_{ji},~\alpha^{k\ell}_{ji},
\alpha_h,c_{ji}\in R,\\
&{~}&\hbox{satisfying~if}~\alpha^{kl}_{ji}\ne 0,~\hbox{then}~i<k\le
\ell <j,~\hbox{and}~k-i=j-\ell .
\end{array}$$
Define two $R$-submodules of the free $R$-module $\RS$:
$$\begin{array}{lll}
{\cal E}_1&=&R\hbox{-Span}\left\{ g_{ji}~\Big |~1\le i<j\le n\right\} ,\\
\\
{\cal E}_2&=&R\hbox{-Span}\left\{
X_ig_{ji},~g_{ji}X_i,~X_jg_{ji},~g_{ji}X_j~\Big |~1\le i<j\le
n\right\} .\end{array}$$ If, for $1\le i<j<k\le n$, every Jacobi sum
$$\begin{array}{rcl} J(X_k,X_j,X_i)&=&\{ X_k,X_j\}X_i-\lambda_{ki}\lambda_{ji}X_i\{ X_k,X_j\}-\\
&{~}&-\lambda_{ji}\{ X_k,X_i\}X_j+\lambda_{kj}X_j\{ X_k,X_i\}+\\
                    &{~}&+\lambda_{kj}\lambda_{ki}\{ X_j,X_i\}X_k-X_k\{
                    X_j,X_i\}\end{array}$$
is contained in ${\cal  E}_1+{\cal E}_2$,  then $A$ is called  a
$q$-{\it enveloping}  algebra. Clearly, enveloping algebras of
$R$-Lie algebras are special $q$-enveloping algebras with $q=1$. In
[Ber], a $q$-PBW theorem for  $q$-enveloping algebras over a
commutative ring was obtained along the line similar to the
classical argument on enveloping algebras of Lie algebras as given
in [Jac], that is, if $A$ is a $q$-enveloping $R$-algebra then $A$
has the PBW $R$-basis $\mathscr{B}=\{ \OV X_1^{\alpha_1}\OV
X_2^{\alpha_2}\cdots \OV X_n^{\alpha_n}~|~\alpha_j\in\NZ\}$.}\par
Now, if we use the $\NZ$-graded monomial ordering
$X_1\prec_{gr}X_2\prec_{gr}\cdots\prec_{gr} X_n$ on $\B_R$ with
respect to the natural $\NZ$-gradation of $\RS$ (i.e., deg$X_i=1$,
$1\le i\le n$), then $\G$ satisfies $\LM (g_{ji})=X_jX_i$ for all
$1\le i<j\le n$. Hence, by Theorem 2.1, the set $\G$ of the defining
relations of a $q$-enveloping $R$-algebra is a monic Gr\"obner basis
for the ideal $I=\langle\G\rangle$ in $\RS$. In particular, all
quantum algebras over $R=\mathbb{C}[[h]]$ which are q-enveloping
algebras appeared in [Ber] are defined by  monic Gr\"obner bases.
{\parindent=0pt\v5

{\bf Remark}  It is necessary to point out that if $R=K$ is a field,
then the fact that the set of defining relations $\G$ of a
$q$-enveloping $K$-algebra $A$ forms a Gr\"obner basis of the ideal
$I=\langle\G\rangle$ was proved in ([Li2], CH.III) directly by using
the termination theorem through the division algorithm. Here our
last example provides the general result for all $q$-enveloping
algebras over an arbitrary commutative ring. \v5

{\bf Example 4.} This example generalizes the previous three
examples but uses an ad hoc monomial ordering. As an application we
show that, over a commutative ring $R$, the PBW generators of the
quantum algebra $U_q^+(A_N)$ derived in [Rin] provides another
Gr\"obner defining set for $U^+_q(A_N)$.}\par With $\RS =R\langle
X_1,...,X_n\rangle$, consider the $R$-algebra $A=\RS
/\langle\G\rangle$ defined by the subset $\G$ consisting of
$\frac{n(n-1)}{2}$ elements
$$\begin{array}{rcl} g_{ji}&=&X_jX_i-q_{ji}X_iX_j-
\sum_{\alpha}\lambda_{\alpha}X_{i_1}^{\alpha_1}X_{i_2}^{\alpha_2}\cdots
X_{i_s}^{\alpha_s}+\lambda_{ji},~1\le i<j\le n,\\
&{~}&\hbox{where}~q_{ji},\lambda_{\alpha},\lambda_{ji}\in
R,~\alpha_k\in\NZ,~i<i_1\le i_2\le\cdots\le i_s<j. \end{array}$$  It
is well-known that numerous iterated skew polynomial algebras over
$R$ are defined subject to such relations, and consequently they
have the PBW $R$-basis $\mathscr{B}=\{ \OV X_1^{\alpha_1}\OV
X_2^{\alpha_2}\cdots \OV X_n^{\alpha_n}~|~\alpha_j\in\NZ\}$. Under
the assumption that $A$ has the PBW $R$-basis as described we aim to
show that $\G$ is a monic Gr\"obner basis of $\langle\G\rangle$. In
view of Theorem 2.1, it is sufficient to introduce a monomial
ordering on $\B_R$ so that $\LM (g_{ji})=X_jX_i$ for all $1\le
i<j\le n$. To this end, let $R[t]=R[t_1,...,t_n]$ be the commutative
polynomial $R$-algebra of $n$ variables. Consider the canonical
algebra epimorphism $\pi$: $\RS\r R[t]$ with $\pi (X_i)=t_i$. If we
fix the lexicographic ordering $X_1<_{lex}X_2<_{lex}\cdots
<_{lex}X_n$ on $\B_R$ of $\RS$ (note that $<_{lex}$ is not a
monomial ordering on $\B_R$) and fix an arbitrarily chosen monomial
ordering $\prec$ on the standard $R$-basis $\mathbb{B}_R=\{
t_1^{\alpha_1}t_2^{\alpha_2}\cdots
t_n^{\alpha_n}~|~\alpha_j\in\NZ\}$ of $R[t]$, respectively, then, as
in [EPS], a monomial ordering $\prec_{et}$ on $\B_R$,  which is
called the {\it lexicographic extension} of the given monomial
ordering $\prec$ on $\mathbb{B}_R$, may be obtained as follows: for
$u,v\in\B_R$,
$$u\prec_{et}v~\hbox{if}~\left\{\begin{array}{l}
\pi (u)\prec\pi (v),\\ \hbox{or}\\
\pi (u)=\pi
(v)~\hbox{and}~u<_{lex}v~\hbox{in}~\B_R.\end{array}\right.$$ In
particular, with respect to the monomial ordering $\prec_{et}$
obtained by using  the lexicographic ordering
$t_n\prec_{lex}t_{n-1}\prec_{lex}\cdots \prec_{lex}t_1$ on
$\mathbb{B}_R$, we see that $\LM (g_{ji})=X_jX_i$ for all $1\le
i<j\le n$, as required by Theorem 2.1.\par In [Rin] it was proved
that $U^+_q(A_N)$ has  $m=\frac{N(N+1)}{2}$ generators $x_1,...,x_m$
satisfying the relations:
$$\begin{array}{rcl} x_jx_i&=&q^{v_{ji}}x_ix_j-r_{ji},~1\le i<j\le m,~
\hbox{where}~v_{ji}=(wt(x_i),wt(x_j)),~\hbox{and}~\\
&{~}&r_{ji}~\hbox{is a linear combination of monomials of the form}~
x_{i+1}^{\alpha_{i+1}}x_{i+2}^{\alpha_{i+2}}\cdots
x_{j-1}^{\alpha_{j-1}},\end{array}$$  and that $U^+_q(A_N)$ is an
iterated skew polynomial algebra generated by $x_1,...,x_m$ subject
to the above relations. Thus $U_q^+(A_N)$ has the PBW basis $\{
x_1^{\alpha_1}x_2^{\alpha_2}\cdots
x_m^{\alpha_m}~|~\alpha_j\in\NZ\}$, and consequently $\G =\{
g_{ji}=x_jx_i-q^{v_{ij}}x_ix_jr_{ji}~|~1\le i<j\le m\}$ forms a
monic Gr\"obner defining set of $U^+_q(A_N)$ with respect to the
monomial ordering $\prec_{et}$ as described before.
{\parindent=0pt\v5

{\bf Remark} If, in the defining relations given in the last
example, the condition $i<i_1\le i_2\le\cdots\le i_s<j$ is replaced
by $1\le i_1\le i_2\le\cdots\le i_s\le i-1$, then a similar result
holds.}\v5

The next three examples provide monic Gr\"obner bases which are not
necessarily the type as described in previous Examples 3 -- 4, but
they all give rise to PBW $R$-bases. {\parindent=0pt\v5

{\bf Example 5.}  Let $R$ be a commutative ring, and let $I$ be the
ideal of the free $R$-algebra $\RS =R\langle X_1,X_2\rangle$
generated by the single element
$$g_{21}=X_2X_1-qX_1X_2-\alpha X_2-f(X_1),$$ where  $q,\alpha\in R$,
and $f(X_1)$ is a polynomial in the variable $X_1$. Assigning to
$X_1$ the degree 1, then in either of the following two cases:\par
(a) deg$f(X_1)\le 2$, and $X_2$ is assigned  the degree 1;\par (b)
deg$f(X_1)=n\ge 3$, and $X_2$ is assigned  the degree n,\par $\G =\{
g_{21}\}$ forms an LM-reduced monic Gr\"obner basis for $I$. For, in
both cases we may use the $\NZ$-graded lexicographic ordering
$X_1\prec_{gr} X_2$ with respect to the natural $\NZ$-gradation of
$\KS$, respectively the weight $\NZ$-gradation of $\RS$ with weight
$\{ 1,n\}$, such that $\LM (g_{21})=X_2X_1$, and then we see that
the only  overlap element of $\G$ is $o(g_{21},1;~1,g_{21})=0$.
Thus, by Theorem 2.1 in both cases the algebra $A=\RS /I$ has the
PBW $R$-basis $\mathscr{B}=\{
\OV{X}_1^{\alpha}\OV{X}_2^{\beta}~|~\alpha ,\beta\in\NZ\}$. \v5

{\bf Example 6.} Let $R$ be a commutative ring, and let $\RS
=R\langle X_1,X_2,X_3\rangle$ be the free $R$-algebra generated by
$X=\{ X_1,X_2,X_3\}$. This example provides a family of algebras
similar to the enveloping algebra $U(\textsf{sl}(2,R))$ of the
$R$-Lie algebra $\textsf{sl}(2,R)$, that is, we consider the algebra
$A=\KS /\langle\G\rangle$ with $\G$ consisting of
$$\begin{array}{l} g_{31}=X_3X_1-\lambda X_1X_3+\gamma X_3,\\
g_{12}=X_1X_2-\lambda X_2X_1+\gamma X_2,\\
g_{32}=X_3X_2-\omega X_2X_3+f(X_1),\end{array}$$ where $\lambda
,\gamma ,\omega\in R$, and $f(X_1)$ is a polynomial in the variable
$X_1$. It is clear that $A=U(\textsf{sl}(2,R))$ in case $\lambda
=\omega =1$, $\gamma =2$ and $f(X_1)=-X_1$.}\par Suppose $f(X_1)$
has degree $n\ge 1$. Then we can always equip $\RS$ with a weight
$\NZ$-gradation by assigning to $X_1$, $X_2$ and $X_3$ the positive
degree $n_1$, $n_2$, $n_3$ respectively (for instance, $(1,1,1)$ if
deg$f(X_1)=n\le 2$; $(1,n,n)$ if deg$f(X_1)=n>2$), such that $\LM
(\G )=\{ X_3X_1,~X_1X_2,~X_3X_2\}$ with respect to the $\NZ$-graded
monomial ordering $X_2\prec_{gr}X_1\prec_{gr}X_3$ on $\B_R$. In the
case that $R=K$ is a field, it was verified in ([Li4], Example 7)
that $\G$ is a Gr\"obner basis for the ideal $\langle\G\rangle$ in
$\KS$ with respect to the same $\prec_{gr}$. Hence, by Proposition
1.6, $\G$ is a Gr\"obner basis for the ideal $\langle\G\rangle$ in
$\RS$. It follows from Theorem 2.1 that the algebra $A=\RS
/\langle\G\rangle$ has the PBW $R$-basis
$\mathscr{B}=\{\OV{X_2}^{\alpha_2}\OV{X_1}^{\alpha_1}\OV{X_3}^{\alpha_3}~
|~\alpha_j\in\NZ\}.$\par Let us point out that in the case that
$f(X_1)$ has degree $\le 2$, i.e., $f(X_1)$ is of the form
$$f(X_1)=aX_1^2+bX_1+c~\hbox{with}~a,b,c\in R,$$
if deg$X_1=$ deg$X_2=$ deg$X_3=1$ is used, the algebra $A$ provides
$R$-versions of some popularly studied algebras over a field $K$,
for instance, {\parindent=0pt
\par
(a) let $\zeta\in R$ be invertible, and put $\lambda =\zeta^4$,
$\omega =\zeta^2$, $\gamma =-(1+\zeta^2)$, $a=0=c$, and $b=-\zeta$,
then $A$ is just the $R$-version of the Woronowicz's deformation of
$U(\textsf{sl}(2,K))$ introduced in the noncommutative differential
calculus;\par (b) if $\lambda\gamma wb\ne 0$ and $c=0$, then $A$ is
just the $R$-version of Le Bruyn's conformal $\textsf{sl}(2,K)$
enveloping algebra [LB] which provides a special family of Witten's
deformation of $U(\textsf{sl}(2,K))$ in quantum group theory.\v5

{\bf Example 7.} Let $\G$ be the subset of the free $R$-algebra $\RS
=K\langle X_1,X_2,X_3\rangle$ consisting of
$$\begin{array}{l} g_{21}=X_2X_1-X_1X_2,\\ g_{31}=X_3X_1-\lambda X_1X_3-\mu X_2X_3-\gamma X_2,\\
g_{32}=X_3X_2-X_2X_3.\end{array}\lambda ,\mu ,\gamma\in R,$$ Then,
under the $\NZ$-graded lexicographic ordering
$X_1\prec_{gr}X_2\prec_{gr}X_3$ with respect to the natural
$\NZ$-gradation of $\RS$, $\LM (g_{ji})=X_jX_i$, $1\le i<j\le 3$,
and the only nontrivial overlap element of $\G$ is
$S_{321}=o(g_{32},X_1;~X_3,g_{21})=-X_2X_3X_1+X_3X_1X_2$. One checks
easily that $\OV{S_{321}}^{\G}=0$. By Theorem 2.1, $\G$ is an
LM-reduced monic Gr\"obner basis for the ideal $\langle\G\rangle$.
Hence, by Theorem 2.1 the algebra $A=\RS /\langle\G\rangle$ has the
PBW $R$-basis
$\mathscr{B}=\{\OV{X_1}^{\alpha_1}\OV{X_2}^{\alpha_2}\OV{X_3}^{\alpha_3}~
|~\alpha_j\in\NZ\}.$}\v5

\section*{3. PBW Isomorphisms and Applications}
In this section we show that the working principle via PBW
isomorphism developed in [LWZ] and [Li3] can be generalized to study
algebras defined by monic Gr\"obner bases over a commutative ring
$R$. All notions and notations used in previous sections are
maintained. \v5

Let $R$ be an arbitrary commutative ring, $\RS =R\langle
X_1,...,X_n\rangle$ the free $R$-algebra of $n$ generators, and
$\B_R$ the standard free $R$-basis of $\RS$. Consider a weight
$\NZ$-gradation of $\RS$ subject to deg$(X_i)=n_i>0$, $1\le i\le n$,
that is, $\RS =\oplus_{p\in\NZ}\RS_p$ with $\RS_p=R$-span$\{
w\in\B~|~ \hbox{deg}(w)=p\}$. For an element $f\in \RS$, say
$f=F_0+F_1+\cdots +F_p$ with $F_i\in \RS_i$ and $F_p\ne 0$, let
$\LH_{\NZ} (F)$ denote the $\NZ$-{\it leading homogeneous element}
of $f$, i.e., $\LH_{\NZ}(f)=F_p$. Then every ideal $I$ of $\RS$ is
associated to an $\NZ$-graded ideal $\langle\LH_{\NZ}(I )\rangle$
generated by the set of $\NZ$-leading homogeneous elements
$\LH_{\NZ}(I )=\{ \LH_{\NZ}(f)~|~f\in I\}$. Adopting the notion and
notation as in [Li3], we call the $\NZ$-graded algebra $A^{\NZ}_{\rm
LH}=\RS /\langle\LH_{\NZ}(I)\rangle$ the $\NZ$-{\it leading
homogeneous algebra} of the algebra $A=\RS /I$. On the other hand,
noticing that $\RS$ is also a $B_R$-graded algebra by the
multiplicative monoid $B_R$, i.e., $\RS =\oplus_{w\in\B_R}\RS_w$
with $\RS_w=Rw$, if $\prec$ is a monomial ordering on $\B_R$ and if
$f=\sum_{i=1}^n\lambda_iw_i\in\RS$ with $w_1\prec
w_2\prec\cdots\prec w_n$, then the term $\lambda_nw_n$ is called the
$\B_R$-{\it leading homogeneous element} of $f$ and is denoted by
$\LH_{\B_R}(f)$. Thus each ideal $I$ of $\RS$ is associated to a
$\B_R$-graded ideal $\langle\LH_{\B_R}(I)\rangle$ generated by the
set of $\B_R$-leading homogeneous elements $\LH_{\B_R}(I)=\{
\LH_{\B_R}(f)~|~f\in I\}$, and similarly, the $\B_R$-graded algebra
$A^{\B_R}_{\rm LH}=\RS /\langle\LH_{\B_R}(I)\rangle$ is referred to
as the $\B_R$-{\it leading homogeneous algebra} of the algebra
$A=\RS /I$. Furthermore, consider the $\NZ$-grading filtration
$F^{\NZ}\RS$ of $\RS$ defined by
$$F^{\NZ}_p\RS =\oplus_{i\le p}\RS_i,\quad p,i\in\NZ ,$$
and the $\B_R$-grading filtration $F^{\B_R}\RS$ of $\RS$ defined by
$$F^{\B_R}_w\RS =\oplus_{u\preceq w}\RS_{u},\quad w,u\in\B_R .$$
If $I$ is an ideal of $\RS$, then the algebra $A=\RS /I$ has the
$\NZ$-filtration $F^{\NZ}A$ induced by $F^{\NZ}\RS$, i.e.,
$$F^{\NZ}_pA=(F^{\NZ}_p\RS +I)/I,\quad p\in\NZ,$$
respectively the $\B_R$-filtration $F^{\B_R}A$ induced by
$F^{\B_R}\RS$, i.e.,
$$F^{\B_R}_wA=(F^{\B_R}_w\RS +I)/I,\quad w\in\B_R.$$
Note that if each $X_i$ has degree 1, $1\le i\le n$, then the
filtration $F^{\NZ}A$ is just the commonly used {\it natural
$\NZ$-filtration}. Let $G^{\NZ}(A)=\oplus_{p\in\NZ}G^{\NZ}(A)_p$
with $G^{\NZ}(A)_p=F^{\NZ}_pA/F^{\NZ}_{p-1}A$ be the associated
$\NZ$-graded algebra of $A$ determined by $F^{\NZ}A$, respectively
$G^{\B_R}(A)=\oplus_{w\in\B_R}G^{\B_R}(A)_w$ with
$G^{\B_R}(A)_w=F^{\B_R}_wA/F^{\B_R}_{\prec w}A$ the associated
$\B_R$-graded algebra of $A$ determined by $F^{\B_R}A$, where
$F^{\B_R}_{\prec w}A=\cup_{u\prec w}F_u^{\B_R}A$. We have the
following analogue of ([Li3], Theorem 1.1). Since the proof of this
result is similar to that given in loc. cit., we omit it here.
{\parindent=0pt\v5

{\bf 3.1. Theorem} With notation as above, there are graded
$R$-algebra isomorphisms:
$$A^{\NZ}_{\rm LH}=\RS /\langle\LH_{\NZ}(I)\rangle\cong G^{\NZ}(A),\quad
A^{\B_R}_{\rm LH}=\RS /\langle\LH_{\B_R}(I)\rangle\cong G^{\B_R}(A)
.$$\par\QED}\v5

Since we are using an arbitrary commutative ring $R$ (instead of a
field) as the coefficient ring,  the  next lemma makes the {\it key
bridge} for us to generalize the working principle of [LWZ] and
[Li3] to quotient algebras of $\RS$ defined by monic Gr\"obner 
bases.   {\parindent=0pt\v5

{\bf 3.2. Lemma} Let $\RS$ be equipped with the fixed weight
$\NZ$-gradation as before, and $I$ an ideal of $\RS$. Put $J=\langle
\LH_{\NZ}(I)\rangle$. The following two statements hold.\par (i) If
$h$ is a nonzero homogeneous element of $\RS$, then $h\in J$ if and
only if $h\in \LH_{\NZ}(I)$. Hence $\LH_{\NZ}(J)=\LH_{\NZ} (I)$.\par
(ii) Let $\prec_{gr}$ be an $\NZ$-graded monomial ordering on $\B_R$
with respect to the fixed weight $\NZ$-gradation of $\RS$.  Then
$\LH_{\B_R}(J)=\LH_{\B_R}(I)$ and $\LM (J)=\LM (I)$.\par (iii) Let
$\prec_{gr}$ be an $\NZ$-graded monomial ordering on $\B_R$ with
respect to the fixed weight $\NZ$-gradation of $\RS$. If $\G$ is a
monic Gr\"obner basis of $I$, then
$$\langle\LH_{\B_R}(J)\rangle =\langle\LH_{\B_R}(I)\rangle =\langle\LM (\G )\rangle
=\langle\LM(I)\rangle =\langle\LM(J)\rangle .$$ \par
{\bf Proof}  (i) Let $h$ be a nonzero homogeneous element in $\RS$.
If $h\in J$, then
$$h=\sum_{i,j} H_{ij}\LH_{\NZ}(f_i)T_{ij},~\hbox{where}~H_{ij},~T_{ij}~\hbox{are homogeneous elements and} ~f_i\in I.$$
If we write $f_i=\LH_{\NZ}(f_i)+f_i'$, where deg$(f_i')<$
deg$(f_i)$, then $f=\sum_{i,j} H_{ij}f_iT_{ij}\in I$ and
$$f=\sum_{ij} H_{ij}\LH_{\NZ}(f_i)T_{ij}+\sum_{i,j} H_{ij}f_i'T_{ij}=h+\sum_{i,j} H_{ij}f_i'T_{ij}.$$
It follows immediately that $h=\LH_{\NZ}(f)\in\LH (I)$. This shows
that  $\LH_{\NZ}(J)\subseteq \LH_{\NZ} (I)$ and hence the equality
holds. \par
(ii) Note that $\prec_{gr}$ is an $\NZ$-graded monomial ordering on
$\B_R$, every element of $\B_R$ is an $\NZ$-homogeneous element, and
thus for $f\in \RS$ we have
$$\LH_{\B_R}(f)=\LH_{\B_R}(\LH_{\NZ}(f))~\hbox{and}~\LM (f)=\LM (\LH_{\NZ}(f))\leqno{(*)}$$
It follows from (i) and the above formula $(*)$ that
$$\begin{array}{c} \LH_{\B_R}(J)=\LH_{\B_R}(\LH_{\NZ}(J))=\LH_{\B_R}(\LH_{\NZ}(I))=\LH_{\B_R}(I),\\
\LM (J)=\LM (\LH_{\B_R}(J))=\LM (\LH_{\B_R}(I))=\LM
(I).\end{array}$$\par
(iii) Let $f\in\RS$ be a monic element with respect to the fixed
monomial ordering $\prec_{gr}$, say $f=w+\sum\lambda_iw_i$ with
$w,w_i\in\B_R$, $\lambda_i\in R$ and $\LM (f)=w$. Then it is clear
that
$$\LM(f)=w=\LH_{\B_R}(f).\leqno{(**)}$$
So, if $\G$ is a monic Gr\"obner basis of $I$ with respect to
$\prec_{gr}$, then the above formula $(**)$ implies $\LM (\G
)=\LH_{\B_R}(\G )\subset\LH_{\B_R}(I)$. Hence, by (ii) and
Proposition 1.2 we obtain the desired equalities:
$$\langle\LH_{\B_R}(J)\rangle =\langle\LH_{\B_R}(I)\rangle =\langle\LM (\G )\rangle
=\langle\LM(I)\rangle =\langle\LM(J)\rangle .$$
 \QED} \v5

Next, we show that an analogue of ([LWZ], Theorem 2.3.2 (i)
$\Leftrightarrow$ (iii)) holds true for  monic Gr\"obner bases in
$\RS$.{\parindent=0pt\v5

{\bf 3.3. Theorem} Let $I$ be an ideal of $\RS$. With notation as
above,  if $\prec_{gr}$ is an $\NZ$-graded monomial ordering on
$\B_R$ with respect to a fixed weight $\NZ$-gradation of $\RS$, the
following two statements are equivalent for a subset $\G\subset
I$:\par (i) $\G$ is a monic Gr\"obner basis of $I$ ;\par (ii)
$\LH_{\NZ}(\G )=\{ \LH_{\NZ}(g)~|~g\in\G\}$ is a monic Gr\"obner
basis for the $\NZ$-graded ideal $\langle\LH_{\NZ}(I
)\rangle$.\vskip 6pt

{\bf Proof} Since we are using the $\NZ$-graded monomial ordering
$\prec_{gr}$ on $\B_R$, by Lemm 3.2 or its proof, a subset $\G$ of
$\RS$ is monic if and only if $\LH_{\NZ}(\G )$ is monic, and we have
$$\langle\LM (I)\rangle =\langle\LM (\G )\rangle~\hbox{if and only if}~
\langle\LM (\langle\LH_{\NZ}(I )\rangle )\rangle =\langle\LM
(\LH_{\NZ}(\G ))\rangle .$$ It follows from Proposition 1.2 that
$\G$ is a monic Gr\"obner basis for the ideal $I$ if and only if
$\LH_{\NZ}(\G )$ is a monic Gr\"obner basis for the $\NZ$-graded
ideal $\langle\LH_{\NZ}(I)\rangle$, proving the equivalence of (i)
and (ii).\QED\v5

{\bf Remark} Also we point out that an analogue of ([LWZ], Theorem
2.3.2 (i) $\Leftrightarrow$ (ii)) works well for monic Gr\"obner
bases when the homogenization of $I$ in $\RS [t]$ is considered.}\v5

Combining the previous 3.1 -- 3.3, we get immediately the result
presenting the associated $\NZ$-graded algebra, respectively the
associated $\B_R$-graded algebra via a monic Gr\"obner
basis.{\parindent=0pt\v5

{\bf 3.4. Theorem} Let $\RS$ be equipped with a fixed weight
$\NZ$-gradation as before, and $I$ an ideal of $\RS$. If $\G$ is a
monic Gr\"obner basis of $I$ with respect to an $\NZ$-graded
monomial ordering $\prec_{gr}$ on $\B_R$, then we have the graded
algebra isomorphisms
$$\begin{array}{l} A^{\NZ}_{\rm LH}=\RS /\langle\LH_{\NZ}(I)\rangle =
\RS /\langle\LH_{\NZ}(\G )\rangle
\cong G^{\NZ}(A),\\
A^{\B_R}_{\rm LH}=\RS /\langle\LH_{\B_R}(I)\rangle =\RS /\langle\LM
(\G
)\rangle \cong G^{\B_R}(A), \\
(A^{\NZ}_{\rm LH})^{\B_R}_{\rm LH}=\RS
/\langle\LH_{\B_R}\langle\LH_{\NZ}(I)\rangle )\rangle =\RS
/\langle\LM (\G )\rangle \cong G^{\B_R}(A^{\NZ}_{\rm
LH}).\end{array}$$\par\QED}\v5

As in [Li3] we call the graded algebra isomorphisms presented in the
last theorem the $\NZ$-{\it PBW isomorphism} (for the first one) and
$\B_R$-{\it PBW isomorphism} (for the last two), determined by the
given monic Gr\"obner basis $\G$ respectively.\v5

Focusing on the first isomorphism of Theorem 3.4, typical examples 
can be given by using the Gr\"obner defining relations of Weyl 
algebras and enveloping algebras of Lie algebras, or more generally, 
the Gr\"obner defining relations of $q$-enveloping algebras 
determined in Example 3 of the last section, over a commutative 
ring. Here we specify several other examples. In all examples given 
below, $R$ is an arbitrary commutative ring.{\parindent=0pt \v5

{\bf Example 1.} Let $X=\{ X_i\}_{i\in J}$ and $\textsf{C}=\RS
/\langle\G\rangle$ the Clifford algebra over $R$, where $\G$
consists of
$$\begin{array}{ll} g_i=X_i^2-q_i,&i\in J,~q_i\in R,\\
g_{k\ell}=X_kX_{\ell}+X_{\ell}X_k-q_{k\ell},&k,\ell\in J,~k>\ell
,~q_{k\ell}\in R.\end{array}$$ Note that if all the $q_i=0 $,
$q_{k,\ell}=0$, we get the defining relations of an $R$-exterior
algebra. It is well known that if $R=K$ is a field, then, under the
$\NZ$-graded lexicographic ordering $\prec_{gr}$ such that
deg$X_i=1$, $i\in J$, and
$$X_{\ell}\prec_{gr}X_k,\quad \ell ,k\in J,~\ell <k,$$
$\G$ forms a Gr\"obner basis for the ideal $\langle\G\rangle$ in
$\KS$ (e.g., see CH.II of [Li2]). It follows from Proposition 1.6
that $\G$ is a Gr\"obner basis for the ideal $\langle\G\rangle$ in
$\RS$. By Theorem 3.4, with respect to the natural $\NZ$-filtration
$F^{\NZ}\textsf{C}$ of \textsf{C}, the associated $\NZ$-graded
algebra $G^{\NZ}(\textsf{C})\cong \RS /\langle\LH_{\NZ}(\G )\rangle$
of \textsf{C} is nothing but an exterior algebra \textsf{E} over
$R$.\v5

{\bf Example 2.} Let $A=R\langle X_1,X_2\rangle /\langle\G\rangle$
be a down-up $R$-algebra in the sense of [Ben], where $\G$ consists
of
$$\begin{array}{l} g_1=X_1^2X_2-\alpha X_1X_2X_1-\beta X_2X_1^2-\gamma X_1,\\
g_2=X_1X_2^2-\alpha X_2X_1X_2-\beta X_2^2X_1-\gamma
X_2,\end{array}\quad \alpha ,\beta \in R.$$  It is well known that
if $R=K$ is a field, then, under the $\NZ$-graded lexicographic
ordering $\prec_{gr}$ such that deg$X_1=$ deg$X_2=1$ and
$X_2\prec_{gr} X_1$, $\G$ forms a Gr\"obner basis for the ideal
$\langle\G\rangle$ in $\KS$ (e.g., see CH.II of [Li2]). It follows
from Proposition 1.6 that $\G$ is a Gr\"obner basis for the ideal
$\langle\G\rangle$ in $\RS$. By Theorem 3.4, with respect to the
natural $\NZ$-filtration $F^{\NZ}A$ of $A$, the associated
$\NZ$-graded algebra $G^{\NZ}(A)\cong R\langle X_1,X_2\rangle
/\langle\LH_{\NZ}(\G )\rangle$ of $A$  is a down-up algebra over $R$
with the set of defining relations $\LH_{\NZ}(\G )=\{\LH_{\NZ}(g_1)
=X_1^2X_2-\alpha X_1X_2X_1-\beta X_2X_1^2, ~ \LH_{\NZ}(g_2)=
X_1X_2^2-\alpha X_2X_1X_2-\beta X_2^2X_1\}$; in particular, one sees
that if $\alpha =2$ and $\beta =-1$, then $G^{\NZ}(A)$ is nothing
but the universal enveloping algebra of the $(-)$-part (or
$(+)$-part) of the Kac-Moody $R$-Lie algebra associated to the
Cartan matrix $\left (\begin{array}{cc} 2&-1\\ -1&2\end{array}\right
)$.\v5

{\bf Example 3.} Let $A=R\langle X_1,X_2\rangle /\langle
g_{21}\rangle$ be the $R$-algebra as given in (Section 2, Example
5). Then by Theorem 3.4, with respect to both the natural
$\NZ$-filtration and the weight $\NZ$-filtration induced by the
weight $\NZ$-grading filtration of $R\langle X_1,X_2\rangle$, $A$
has the associated $\NZ$-graded algebra $G^{\NZ}(A)\cong R\langle
X_1,X_2\rangle /\langle X_2X_1-qX_1X_2\rangle$, which, in the case
that $q$ is invertible, is the coordinate ring of the quantum plane
over $R$. \v5

{\bf Example 4.} Let $A=R\langle X_1,X_2, X_3\rangle /\langle
\G\rangle$ be the $R$-algebra as given in (Section 2, Example 6). In
the case that $f(X_1)$ has degree $\le 2$, i.e., $f(X_1)$ is of the
form
$$f(X_1)=aX_1^2+bX_1+c~\hbox{with}~a,b,c\in R,$$
then by Theorem 3.4, with respect to the natural $\NZ$-filtration
$F^{\NZ}A$, $A$ has the associated $\NZ$-graded algebra
$G^{\NZ}(A)\cong R\langle X_1,X_2,X_3\rangle /\langle\LH_{\NZ}(\G
)\rangle$ with
$$\LH_{\NZ} (\G )=\{X_3X_1-\lambda X_1X_3,~X_1X_2-\lambda X_2X_1,~X_3X_2-\omega X_2X_3+aX_1^2\} ;$$
while in the case that $f(X_1)$ has degree $n\ge 3$, if the weight
$(1,n,n)$ is used, then by Theorem 3.4, with respect to the  weight
$\NZ$-filtration $F^{\NZ}A$ induced by the weight $\NZ$-grading
filtration of $R\langle X_1,X_2\rangle$, $A$ has the associated
$\NZ$-graded algebra $G^{\NZ}(A)\cong R\langle X_1,X_2,X_3\rangle
/\langle\LH_{\NZ}(\G )\rangle$ with
$$\LH_{\NZ} (\G )=\{X_3X_1-\lambda X_1X_3,~X_1X_2-\lambda X_2X_1,~X_3X_2-\omega X_2X_3\}$$
}\v5

By referring to the well-known filtered-graded comparison principle 
for algebras with an $\NZ$-filtration ([MR], [Li1], [LVO1], [Li3]), 
we now summarize, without proof, several  applications of Theorem 
3.3 and Theorem 3.4. Let $R$ be an arbitrary commutative ring. For 
the convenience, in what follows we let the free $R$-algebra 
$\RS=R\langle X_1,...,X_n\rangle$ be equipped with a fixed weight 
$\NZ$-gradation, $I=\langle\G \rangle$ an ideal of $\RS$  generated 
by a monic Gr\"obner basis $\G$ with respect to an $\NZ$-graded 
monomial ordering $\prec_{gr}$ on the standard $R$-basis $\B_R$ of 
$\RS$, and $A=\RS /I$. Then the following diagram may indicate how 
all results to be given will work:
$$\begin{diagram} &&A=\RS /\langle\G\rangle&&\\
&\NE^{\textsf{lifting}}&&\NW^{\textsf{lifitng}}&\\
G^{\B_R}(A)&&&&G^{\NZ}(A)\mapleft{\cong}{\scriptstyle{\NZ {\small\rm
-PBW}}}
\FRAC{\RS}{\langle\LH_{\NZ}(\G )\rangle}=A^{\NZ}_{\rm LH}\\
\hbox{~~~~~~~~}\uTo^{\cong}_{\scriptstyle{\B_R {\rm\small -PBW}}}&&&&~~\uTo_{\textsf{~lifting}}\\
A^{\B_R}_{\rm LH}=\FRAC{\RS}{\langle\LM (\G )\rangle}&&
\rTo^{~~~~\cong}_{\scriptstyle{~~~~~~\B_R {\rm\small
-PBW}}}&&G^{\B_R}(A^{\NZ}_{\rm LH})
\end{diagram}$${\parindent=0pt\par

{\bf 3.5. Theorem}  Under the  respective canonical algebra
epimorphism, the set $N(\G )$ of normal monomials in $\B_R$ (mod
$\G$), projects to a free $R$-basis for the algebras $A=\RS /I$,
$A^{\NZ}_{\rm LH}=\RS /\langle\LH_{\NZ} (I)\rangle$ , and
$A^{\B_R}_{\rm LH}=\RS /\langle\LM (I)\rangle$ respectively, and
thereby to a free $R$-basis for $G^{\NZ}(A)$, $G^{\B_R}(A)$, and
$G^{\B_R}(A^{\NZ}_{\rm LH}(A))$, respectively.\QED  \v5

{\bf 3.6. Theorem} Bearing $A^{\B_R}_{\rm LH}=\RS /\langle\LM (\G
)\rangle$ in mind, the following statements hold.\par (i) If
$A^{\B_R}_{\rm LH}$ is a (semi-)prime ring, then $A^{\NZ}_{\rm LH}$
is a (semi-)prime ring (hence $G^{\NZ}(A)$ is a (semi-)prime ring),
and $A$ is a (semi-)prime ring.\par (ii) If $A^{\B_R}_{\rm LH}$ is
$\B_R$-graded left Noetherian, that is, every $\B_R$-graded left
ideal of $G(A)$ is finitely generated, then $A^{\NZ}_{\rm LH}$ is
left Noetherian (hence $G^{\NZ}(A)$ is left Noetherian), and $A$ is
left Noetherian.\par (iii) If $A^{\B_R}_{\rm LH}$ is $\B_R$-graded
left Artinian, that is, $A^{\B_R}_{\rm LH}$ satisfies the descending
chain condition for $\B_R$-graded left ideals, then $A^{\NZ}_{\rm
LH}$ is left Artinian (hence $G^{\NZ}(A)$ is left Artinian), and $A$
is left Artinian.\par (iv) If $A^{\B_R}_{\rm LH}$ is a $\B_R$-graded
simple $R$-algebra, that is, $A^{\B_R}_{\rm LH}$ does not have
nontrivial $\B_R$-graded ideal, then $A^{\NZ}_{\rm LH}$ is a simple
$R$-algebra (hence $G^{\NZ}(A)$ is a simple $R$-algebra), and $A$ is
a simple $R$-algebra.\par (v) If the Krull dimension (K.dim in the
sense of Gabriel and Rentschler, e.g. see [MR] for the definition)
of $A^{\B_R}_{\rm LH}$ is well-defined, then the Krull dimension of
$A^{\NZ}_{\rm LH}$ (hence of $G^{\NZ}(A)$) and $A$ is defined and
K.dim$A\le$ K.dim $G^{\NZ}(A)=$ K.dim$A^{\NZ}_{\rm LH}\le$
K.dim$A^{\B_R}_{\rm LH}$.
\par
(vi) If $A^{\B_R}_{\rm LH}$ is semisimple (simple) Artinian, then
$A^{\NZ}_{\rm LH}$ is semisimple (simple) Artinian (hence
$G^{\NZ}(A)$ is semisimple (simple) Artinian), and $A$ is semisimple
(simple) Artinian.\par
(vii) Let gl.dim abbreviate the phrase ``global homological
dimension". We have gl.dim$A\le$ gl.dim$G^{\NZ}(A)=$
gl.dim$A^{\NZ}_{\rm LH}\le$  gl.dim$A^{\B_R}_{\rm LH}$.\par
(viii) If $A^{\B_R}_{\rm LH}$ is left hereditary, then $A^{\NZ}_{\rm
LH}$ is  left hereditary (hence $G^{\NZ}(A)$ is left hereditary),
and $A$ is left hereditary.\par
(ix) Let gl.wdim abbreviate the phrase ``global week homological
dimension". We have $\hbox{gl.wdim}A\le$ gl.wdim$G^{\NZ}(A)=$
gl.wdim$A^{\NZ}_{\rm LH}\le~ \hbox{gl.wdim}A^{\B_R}_{\rm LH}$.\par
(x) If $A^{\B_R}_{\rm LH}$ is a Von Neuman regular ring, then
$A^{\NZ}_{\rm LH}$ is  Von Neuman regular ring (hence $G^{\NZ}(A)$
is a Von Neuman regular ring), and  $A$ is a Von Neuman regular
ring.\QED \v5

{\bf 3.7. Theorem} Bearing $A^{\NZ}_{\rm LH}=\RS
/\langle\LH_{\NZ}(\G )\rangle$ in mind, if the role of
$A^{\B_R}_{\rm LH}$ is replaced by $A^{\NZ}_{\rm LH}$, then the
analogues of Theorem 4.6 (i) -- (x) hold true. Moreover, we
have:\par (i) If $A^{\NZ}_{\rm LH}$ is a domain, then $A$ is a
domain.\par (ii) If $A^{\NZ}_{\rm LH}$ is a Noetherian domain and
maximal order in its quotient ring (see e.g. [MR] for the
definition), then $A$ is a Noetherian domain and maximal order in
its quotient ring.\par (iii) If $A^{\NZ}_{\rm LH}$ is an Auslander
regular ring (see e.g. [Li1], [LVO] for the definition), then $A$ is
an Auslander regular ring.\QED}\v5

\v5 \centerline{References}\par \parindent=1truecm\par
\item{[An]} D. J. Anick, On the homology of associative algebras,
{\it Trans. Amer. Math. Soc}., 2(296)(1986), 641--659.
\item{[Ben]} G.~Benkart, Down-up algebras and Witten's deformations of
the universal enveloping algebra of $sl_2$, {\it Contemp. Math.},
224(1999), 29--45.
\re{[Ber]} R. Berger, The quantum Poincar\'e-Birkhoff-Witt theorem,
{\it Comm. Math. Physics}, 143(1992), 215--234.
\item{[Ber1]} G. Bergman, The diamond lemma for ring theory, {\it
Adv. Math}., 29(1978), 178--218.
\item{[BM]} L. A. Bokut and P. Malcolmson, Gr¡§obner-Shirshov basis for
quantum enveloping algebras, {\it Israel Journal of Mathematics},
96(1996), 97--113.
\item{[CE]} B.L. Cox and T.J. Enright, Representations of quantum groups
defined over commutative rings, {\it Comm. Alg}, 23(6)(1995), 2215 -
2254.
\item{[CSS]} Y. Chen, H. Shao and K. P. Shum, Rosso-Yamane Theorem on PBW basis of
$U_q(A_N)$, {\it CUBO, A Mathematical Journal}, 10(3)(2008),
171--194. arXiv:0804.0954v1 [math.RA]

\item{[CU]} S. Cojocaru and V. Ufnarofski, BERGMAN under MS-DOS and
Anick's resolution, {\it Discrete Mathematics and Theoretical
Computer Science}, 1(1997), 139--147.

 \re{[EPS]} D. Eisenbud, I. Peeva and B.
Sturmfels, Non-commutative Gr\"obner bases for commutative algebras,
{\it Proc. Amer. Math. Soc.}, 126(1998), 687--691.

\item{[GI-L]} T.~Gateva-Ivanova and V.~Latyshev, On recognizable
properties of associative algebras,  {\it  J. Symbolic Computation},
6(1988), 371--388.

\item{[Gol]} E. S. Golod, On noncommutative Gr\"obner bases over rings,
 {\it Journal of Mathematical Science}, 140(2)(2007), 239--242.
\item{[Gr]} E. L. Green, Noncommutative Gr¡§obner bases and projective
resolutions, in: {\it Proceedings of the Euroconference
Computational Methods for Representations of Groups and Algebras},
Essen, 1997, (Michler, Schneider, eds), Progress in Mathematics,
Vol. 173, Basel, Birkha¡§user Verlag, 1999, 29--60.
\item{[Jac]} N. Jacobson, {\it Lie Algebras}, New York, Wiley, 1962.
\re{[LB]} L.~Le Bruyn, Conformal $sl_2$ enveloping algebras, {\it
Comm. Alg.}, 23(1995), 1325--1362.
\item{[Li1]} H. Li, {\it Noncommutative Zariskian Filtered Rings}, PhD Thesis, Antwerp University, 1989.
\item{[Li2]} H. Li, {\it Noncommutative Gr\"obner Bases and
Filtered-Graded Transfer}, LNM, 1795, Springer-Verlag, 2002.
\item{[Li3]} H. Li, $\Gamma$-leading homogeneous algebras and Gr\"obner bases,
in: {\it Advanced Lectures in Mathematics}, Vol.8, International
Press, Boston, 2009, 155--200.
\item{[Li4]} H. Li, Looking for Grobner basis theory for (almost) skew 2-nomial
algebras, {\it Journal of Symbolic Computation}, 45(2010), in press, 
online 10 may 2010, doi:10.1016/j.jsc.2010.05.002. (also see 
arXiv:math.RA/0808.1477)
\item{[LVO1]} H. Li and F. Van Oystaeyen, {\it Zariskian Filtrations}, $K$-Monograph in
Math., Vol.2, Kluwer Academic Publishers, 1996.
\item{[LVO2]} H. Li and F. Van Oystaeyen, Reductions and global dimension of quantized algebras over a
regular commutative domain, {\it Comm. Alg.}, 26(4)(1998),
1117--1124.
\item{[LWZ]} H. Li, Y. Wu and J. Zhang, Two applications of
noncommutative Gr\"obner bases, {\it Ann. Univ. Ferrara - Sez. VII -
Sc. Mat.}, XLV(1999), 1--24.
\re{[Man]} Yu.I.~Manin, {\it Quantum Groups and Noncommutative
Geometry}, Les Publ. du Centre de R\'echerches Math., Universite de
Montreal, 1988.
\item{[MR]} J.C.~McConnell and J.C.~Robson, {\it Noncommutative
Noetherian Rings}, John Wiley \& Sons, 1987.
\item{[Mor]} T. Mora, An introduction to commutative and
noncommutative Gr\"obner bases, {\it Theoretic Computer Science},
134(1994), 131--173.
\item{[Rin]} C. M. Ringel, PBW-bases of quantum groups, {\it J. Reine Angew}. Math.
470 (1996), 51--88.
\item{[Ros]} M. Rosso, An analogue of the Poincare-Birkhoff-Witt
theorem and the universal R-matrix of $U_q(sl(N + 1))$, {\it Comm.
Math. Phys}., 124(2)(1989), 307--318.

\item{[Uf1]} V. Ufnarovski, A growth criterion for graphs and
algebras defined by words, {\it Mat. Zametki}, 31(1982), 465--472
(in Russian); English translation: {\it Math. Notes}, 37(1982),
238--241. \item{[Uf2]} V. Ufnarovski, On the use of graphs for
computing a basis, growth and Hilbert series of associative
algebras, (in Russian 1989), {\it Math. USSR Sbornik},
11(180)(1989), 417-428.

\item{[Yam]} I. Yamane, A Poincare-Birkhoff-Witt theorem for quantized
universal enveloping algebras of type $A_N$, Publ., {\it RIMS. Kyoto
Univ}., 25(3)(1989), 503--520.

\end{document}